\documentclass[conference]{IEEEtran}
\IEEEoverridecommandlockouts                              

\usepackage{graphicx}
\usepackage{subfigure}
\usepackage{amssymb}
\usepackage{widetext}
\usepackage{multicol}
\usepackage{amsmath}
\usepackage{lipsum}  
\usepackage{cite}
\usepackage{textcomp}
\def\BibTeX{{\rm B\kern-.05em{\sc i\kern-.025em b}\kern-.08em
    T\kern-.1667em\lower.7ex\hbox{E}\kern-.125emX}}
    
\usepackage{color}

\providecommand{\com[2]}{\begin{tt}[#1: #2]\end{tt}}

\newcommand{\bydef}{\stackrel{\Delta}{=}}

\newcommand{\beq}{\begin{equation}}
\newcommand{\eeq}{\end{equation}}
\newcommand{\beqa}{\begin{eqnarray}}
\newcommand{\eeqa}{\end{eqnarray}}
\newcommand{\beqan}{\begin{eqnarray*}}
\newcommand{\eeqan}{\end{eqnarray*}}
\newcommand{\bef}{\begin{figure}}
\newcommand{\enf}{\end{figure}}

\newcommand{\bi}{\begin{itemize}}
\newcommand{\ei}{\end{itemize}}
\newcommand{\bc}{\begin{center}}
\newcommand{\ec}{\end{center}}
\newcommand{\ba}{\begin{array}}
\newcommand{\ea}{\end{array}}
\newcommand{\be}{\begin{equation}}
\newcommand{\ee}{\end{equation}}
\newcommand{\beno}{\begin{equation*}}
\newcommand{\eeno}{\end{equation*}}
\newcommand{\beqna}{\begin{eqnarray}}
\newcommand{\eeqna}{\end{eqnarray}}
\newcommand{\bd}{\begin{displaymath}}
\newcommand{\ed}{\end{displaymath}}
\newcommand{\beqnd}{\begin{eqnarray*}}
\newcommand{\eeqnd}{\end{eqnarray*}}

\renewcommand{\ni}{\noindent}

\newcommand{\Rank}{\mathrm{Rank}}

\newcommand{\cqfd}{\hfill \rule{2mm}{2mm}\medbreak\indent}

\providecommand{\quotes[1]}{``#1"}

\providecommand{\prt}[1]{\left( #1 \right)}

\newtheorem{theorem}{\bf Theorem}[section]
\newtheorem{lemma}{\bf Lemma}[section]

\newtheorem{proposition}{\bf Proposition}[section]
\newtheorem{corollary}{\bf Corollary}[section]

\newtheorem{definition}{\bf Definition}

\usepackage{color}

\definecolor{red}{rgb}{1,0,0}

\definecolor{blu}{rgb}{0,0,1}

\definecolor{gre}{rgb}{0,0.7,0.3}

\include{temp_jh_ref.bib}

\begin{document}

 \title{Identifiability of dynamical networks with partial node measurements*\\
\thanks{This work is supported
 by the Program Science Without Borders, CNPq - Conselho Nacional de Desenvolvimento Cient'\'{\i}fico e Tecnol\'{o}—gico, Brazil,  by the Belgian Programme on
Interuniversity Attraction Poles, initiated by the Belgian Federal
Science Policy Office, by Wallonie-Bruxelles International, and by a Concerted Research Action (ARC) of the
French Community of Belgium.}
}

\author{\IEEEauthorblockN{1\textsuperscript{st} Julien M. Hendrickx}
\IEEEauthorblockA{\textit{ICTEAM}\\
\textit{ Universit\'e catholique de Louvain  }\\
B1348 Louvain la Neuve, Belgium \\
Julien.Hendrickx@uclouvain.be}
\and
\IEEEauthorblockN{2\textsuperscript{nd} Michel Gevers}
\IEEEauthorblockA{\textit{ICTEAM}\\
\textit{ Universit\'e catholique de Louvain  }\\
B1348 Louvain la Neuve, Belgium \\
Michel.Gevers@uclouvain.be}
\and
\IEEEauthorblockN{3\textsuperscript{rd} Alexandre S.  Bazanella}
\IEEEauthorblockA{\textit{Department of Automation and Energy}\\ 
\textit{ Universidade Federal do Rio Grande do Sul} \\
Porto Alegre, Brazil \\
bazanella@ufrgs.br}
}

\maketitle

\begin{abstract}
Much recent research has dealt with the identifiability of a dynamical network in which the node signals are connected by causal linear  transfer functions and are  excited by known external excitation signals and/or unknown noise signals. A major research  question concerns the identifiability of the whole network - topology and all transfer functions - from the measured node signals and external excitation signals.  So far all  results on this topic have assumed that all node signals are measured.  This paper  presents  the first results for the situation where not all node signals are measurable, under the assumptions that (1) the topology of the network is known, and (2) each node is excited by a known external excitation.  Using graph theoretical properties, we show that the transfer functions that can be identified depend essentially on the topology of the paths linking the corresponding vertices to the measured nodes. A practical outcome  is that, under those assumptions,  a network can often be identified using only a small subset of node measurements. 
\end{abstract}

\begin{IEEEkeywords}
Network Analysis and Control; System identification.
 \end{IEEEkeywords}

\section{Introduction}\label{intro}

This paper examines the identifiability of dynamical networks in which the node signals  are connected by causal linear time-invariant transfer functions and are  excited by known external excitation signals.  Such networks can be looked upon as connected directed graphs in which the edges between the nodes (or vertices) are scalar  transfer functions, and in which known external excitation signals enter into the nodes.

The identification of  networks of linear time-invariant dynamical systems based on the measurement of all   its node signals and of all known external excitation signals acting on the nodes has been the subject of much recent research \cite{Goncalves&Warnick:08,Materassi&Innocenti:10,Dankers&Vandenhof&Heuberger&Bombois:12,Chiuso&Pillonetto:12, Weerts&Dankers&Vandenhof:15,Hayden&Chang&Goncalves&Tomlin:16,Gevers&Bazanella&Parraga:17}. It has been shown in \cite{Goncalves&Warnick:08,Weerts&Dankers&Vandenhof:15,Gevers&Bazanella&Parraga:17} that  identifiability can only be obtained provided prior knowledge is available about the structure of the network, and in particular the structure of the excitation. It is often the case that the  excitation structure is known, i.e. one  often knows at which nodes  external excitation signals  are applied.  A number of conditions for the identifiability of the whole network have been derived under prior assumptions on the structure of the network, involving either its external excitation structure,  or possibly also its internal structure \cite{Goncalves&Warnick:08,Weerts&Dankers&Vandenhof:15,Hayden&Chang&Goncalves&Tomlin:16,Gevers&Bazanella&Parraga:17}.

In all the results accumulated so far on the identifiability of  a  network of dynamical systems, it is assumed that all node signals  are measured.  In this paper we examine the situation where not all node signals are measured, but where the topology of the network is known; this means that the user knows a priori which nodes are connected by nonzero transfer functions.  We also make the simplifying assumption that at each node a known external excitation is applied. In this context, a number of questions can be raised, such as
\begin{enumerate}
\item Can one identify the whole network with a restricted number of node measurements?
\item If so,  are there a minimal  number of nodes that need to be measured?
\item Are some nodes indispensable, in the sense that it is impossible to identify the network without measuring these nodes?
\item If one wants to identify a specific transfer function, can the topology tell us which node or nodes need to be measured?
\item Which transfer functions can be identified from the measure of a  specific subset of nodes?
\end{enumerate}

To answer these questions we shall heavily rely on properties from graph theory, using the connected directed graph corresponding to our network as our major tool. 

To the best of our knowledge, the only other contributions that consider identification in networks using only a subset of measured nodes are \cite{Materassi&Salapaka:15}, \cite{Dankers&Vandenhof&Bombois&Heuberger:16}, \cite{Linder&Enqvist:17}.
However, the problem treated in these papers consists of 
the identification of a subset of the network's transfer functions - typically a single one - and hence
 is only one of the subproblems presented in this paper.
In \cite{Materassi&Salapaka:15} networks driven only by a vector of white noises are considered, i.e. no known external excitation is available. Using the notion of $d$-separation of graphs, the authors derive sufficient conditions on which node signals need to be observed in order to guarantee the identifiability of a desired transfer function link. In \cite{Dankers&Vandenhof&Bombois&Heuberger:16} the objective is also to identify a specific transfer function link (or module), in a network that does have both known external excitation signals and/or noise signals on the nodes. The authors present sufficient conditions  for the selection of a set of measured node signals that will lead to the consistent identification of the desired module. The approach taken in \cite{Linder&Enqvist:17} is quite different. It consists of estimating the desired but unobservable nodes from nodes that are measurable and contain information about them.
Sparse measurements have also been considered in a different context in \cite{Mauroy&Hendrickx:16};
the goal there was to recover the network structure under the assumptions that the local dynamics are known,
as opposed to re-identifying the dynamics and/or the whole network structure.

Our main contribution is to provide necessary and sufficient conditions under which all transfer functions of the network, or a subset of transfer functions, or a single transfer function can be identified from a given set of measured nodes, under the  standing assumptions that the topology of the network is known  and that  there is a known external excitation on each node. Our results are existence results about identifiability; they are not algorithms for the estimation of the transfer functions.  They all take  the form of conditions on the topology of the graph associated to the network. We  also present the computational complexity that is required to check these necessary and sufficient conditions.

In Section~\ref{probstatement} we  first describe the standard network  matrix identifiability problem where all nodes are measured but where the topology is unknown and needs to be identified from data. We   explain that  without any knowledge of the topology the identification of the network's transfer functions from partial node measurements has no solution.  We then show that,  in order to relate the identifiability of a set of transfer functions to the selection of a set of measured nodes on the basis of the network topolgy, one needs to introduce the notion of generic identifiability.  This notion is described intuitively in Section~\ref{probstatement} together with a motivating example.

We  then motivate the reason for addressing the  problem  of network identifiability with partial node measurements in  Section~\ref{motivation} by analyzing three different 3-node networks. We  show  that the nodes that need to be measured to identify all transfer functions depend on the topology of the network and that, in some cases, a unique measurement suffices to identify the whole network. This already yields a positive answer to question 1 above. Our brief analysis of 3-node networks  then leads us, in Section~\ref{basicresults}, to formulate a number of basic results pertaining to questions 2, 3 and 4 above.
We also provide identifiability results for networks that have a special structure, such as a  tree or a loop.

In Section~\ref{pathbased} we  focus on the identifiability of the transfer functions leaving a specific node $i$, i.e. the  transfer functions $G_{ji}$ that connect node $i$ to its outgoing nodes. Our main result in that Section is a necessary and sufficient condition for the identifiability of a  set of transfer functions leaving node $i$. This set of transfer functions is shown to be identifiable from a given set of measured nodes if and only if there are disjoint paths going from these outgoing nodes of $i$ to the set of measured nodes. 

In Section~\ref{measurementbased} we  address question 5 above. Instead of looking at a specific node within the network and examining its paths to a measured node or a set of measured nodes, as was done in Section~\ref{pathbased}, we consider the converse approach. We consider a specific set of measured nodes and we ask which transfer functions can be identified from it. Our main result is a necessary and sufficient condition under which the whole network can be identified from a given set of node measurements. 

In Section~\ref{algocom} we examine the computational complexity of the algorithm to check the identifiability of the whole network or parts of it from a given set of measurements. We show for example that checking the identifiability of the whole network can be achieved at a computational cost of the order of  $L^2 \times n$ where $L$ is the number of nodes and $n$ the number of unknown transfer functions in the network.

In Section~\ref{conclusion} we will conclude and describe  some challenging  open problems that remain to be solved.\\

\section{Statement of the problem}\label{probstatement}
The problem studied in this paper is part of the recent research on the question of identifiability of networks of dynamical systems. We first present the network structure and explain the network identifiability problem as it has so far been posed, i.e.  with all nodes measured. We then pose a new network identifiability problem for the case when not all nodes are
measured. 

We adopt the standard network structure of  \cite{Weerts&Dankers&Vandenhof:15,Gevers&Bazanella&Parraga:17}
for networks whose edges are labeled with scalar proper transfer functions. 
Thus, we consider that the network is made up of $L$ nodes, with  node signals  denoted $\{w_1(t), \ldots, w_L(t)\}$,
and that  these node signals are related to each other and to  external excitation signals $r_j(t), j=1,\dots,L$ by the following
network equations, which we call the {\bf network model} and in which the matrix $G^0(q)$ is called the {\bf  network matrix}:
\be
w(t) = G^0(q) w(t) +  r(t) +v(t) . \label{netmodel} 
\ee
In (\ref{netmodel}) $q^{-1}$ is the delay operator,
$w(t) = [w_1(t), \ldots, w_L(t)]^T$ is the vector of node signals, $r(t) = [r_1(t), \ldots, r_L(t)]^T$ is a vector of known external excitation signals,  $v(t) = [v_1(t), \ldots, v_L(t)]^T$ is 
a vector of stochastic processes, and the dynamic network matrix is of the form
\beqnd 
G^0(q) =\! \left[\! \begin{array}{cccc}0 & \!\! G_{12}^0(q) & \!\! \ldots &\!\!  G_{1L}^0(q) \\G_{21}^0(q) & \!\! 0  \!\!  & \ddots & \!\! G_{2L}^0(q) \\Ê\vdots &\!\! \ddots &\!\!  \ddots &\!\!  \vdots\\
G_{L1}^0(q) & \!\! G_{L2}^0 (q)& \!\! \ldots & \!\! 0\end{array}\! \right] 
\eeqnd

The network (\ref{netmodel}) is assumed to have the following  properties.
\begin{itemize}
\item $G_{ij}^0(q)$ are proper rational transfer functions 
\item the network is well-posed, that is $(I-G^0(q))^{-1}$ is proper and stable  \cite{Araki&Saeki:83}
\item there is a known  external excitation signal $r_i(t)$ on each node; these are available to the user in order to produce informative experiments for identification 
\item the network is weakly connected\footnote{A precise definition will be given in Section~\ref{basicresults}.}
\end{itemize}
In most papers on identifiability of networks based on measurements of all the nodes, the vector $r(t)$ of external excitation signals traditionally enters the nodes via a transfer function matrix $K^0(q)$, i.e. the driving term is $r(t) = K^0(q)\tilde r(t)$ where $\tilde r(t)$ is the vector of external excitations. In this paper on network identifiability using partial node measurements, we adopt the simplified network model (\ref{netmodel}) where $K^0(q)=I$. Observe that, by a simple change of variables, this is equivalent to assuming that   in the traditional model the excitation matrix $K^0(q)$ is known and of full rank. The reason for making this simplifying assumption is that, as we shall see, the problem treated in this paper, even with this assumption, is complex enough and reveals significant  new insights. We expect to be able  to relax this assumption in future work.

The network model (\ref{netmodel}) can be rewritten in a more traditional input-output (I/O) form as follows:
 \be 
w(t) = T^0(q) r(t)  + \bar v(t)\label{iomodel2}
\ee
where 
\begin{eqnarray} 
&& T^0(q) \bydef (I - G^0(q))^{-1}\label{Tdef} \\
&& \bar v(t) = (I - G^0(q))^{-1}v(t). \label{barvdef}
\end{eqnarray}

In this paper we address the question of the identifiability of the network matrix $G^0(q)$
 for the case where not all nodes are measured,
but where the topology of the network is known. 
The reason for the assumption on known topology is that, as  we shall show in Theorem~\ref{Lminus1}, when not all nodes are measured, some knowledge of the topology is required in order to identify the whole network  (in the absence thus of any prior knowledge on the specific transfer functions $G_{ij}(q)$).

Thus we assume that we know that certain transfer functions $G_{ij}(q)$ are zero, and we say that a network matrix $G(q)$ is \emph{consistent} with the topology if it satisfies these constraints. Moreover, we consider that,  
 together with the network (\ref{netmodel}), there is a measurement equation
\be \label{measures}
y(t) = C w(t)
\ee
where $C$ is a $pÊ\times L$ matrix that reflects the selection of measured nodes. That is,
each row of $C$ contains one element $1$ and $L-1$ elements $0$.
We shall denote by $\cal C$ the corresponding subset of nodes selected by $C$.

In this setting, the network under study is given by
\beqna \label{netmodel1}
w(t) & = & G^0(q) w(t) +  r(t) + v(t) \label{netmodel3} \\
y(t) &=&  Cw(t)  \label{netmodel2}
\eeqna
which, in the input-output form, becomes
\be\label{netmodel4}
y(t) = CT^0(q) r(t) + C\bar v(t)
\ee
with $Tï(q)$ and  $\bar v(t)$ defined by (\ref{Tdef}) and (\ref{barvdef}).

We now describe the {\bf network matrix identifiability} problem for such networks; we start by summarizing
the  assumptions that are made throughout this paper.

\ni {\bf Standing assumptions.}
\begin{itemize}
\item The networks we examine are described by (\ref{netmodel3})-(\ref{netmodel2}).
\item The network matrix $G^0(q)$ has the properties defined above and its topology is
known, i.e. one knows {\em a priori} that some of the $G_{ij}^0(q)$ are zero.
\item The excitation vector $r(t)$ is sufficiently rich such 
that   $CT^0(q)$ can be consistently estimated by standard identification
of the open loop MIMO I/O model (\ref{netmodel4}). 
\end{itemize}

Since, for a given $C$, the matrix $CT^0(q)$ can be consistently identified from $\{y(t), r(t)\}$ data, it will
be assumed to be known exactly. The {\bf network matrix identifiability} problem is  whether or not, under the standing assumptions,  one can uniquely recover $G^0(q)$ from $CT^0(q)$.
Specific questions related to network identifiability that are  addressed in this paper are then: 
\begin{itemize}
\item for a given $C$, which transfer functions $G^0_{ij}(q)$ can be uniquely recovered from $CT^0(q)$? 
\item under what conditions can we identify the whole network matrix $G^0(q)$ from $CT^0(q)$?
\end{itemize}

The identification of the transfer functions $G_{ij}(q)$ from $CT^0(q)$ rests on the following relationship
\be \label{TzeroT1}
C T^0(q) = CT(q) = C(I-G(q))^{-1}
\ee
or, equivalently,
\be \label{TzeroT2}
C T^0(q) (I-G(q)) = C
\ee
Since $CT^0(q)$ is assumed known, the question is whether the desired $G_{ij}(q)$ can be uniquely obtained by
solving (\ref{TzeroT2}) for these unknowns, using the knowledge of the network topology. More precisely, we say that the network matrix $G^0$ is identifiable from the measurements $C$ if it is the unique solution of \eqref{TzeroT2} consistent with the topology. Similarly, a specific transfer function in $G^0$ is identifiable from the measurements $C$ if $G_{ij}=G^0_{ij}$ for any solution $G$ of \eqref{TzeroT2} consistent with the topology.

Deciding whether $G^0(q)$ is uniquely recoverable from the identified and exact $CT^0(q)$ can thus be done by checking whether the solution $G(q)$ of (\ref{TzeroT2}) is unique. However, this is of limited interest because it does not take account of the information we have about the known topology of the network. Our ambition in this paper is to make statements about the identifiability of $G^0(q)$ for a given selection $\cal C$ of node measures before we actually compute $CT^0(q)$ from data or even collect the data, i.e. statements that are based  not on the actual numbers that appear in the transfer functions of $CT^0(q)$, but on the topology of the network that is assumed to be known, and which can be represented by a graph associated to $G^0(q)$. 
 
As a consequence, we will introduce the notion of {\it generic identifiability} of $G^0(q)$ because the topology tells us which of the $G_{ij}(q)$ can be nonzero, which impacts on the generic rank of $CT^0(q)$  and of its submatrices\footnote{By the generic rank of a submatrix of $CT^0(q)$ we mean its rank for almost all $G^0(q)$ that are consistent with its associated graph.}, but one cannot exclude the possible situation where a given $G^0(q)$, that is consistent with the topology, happens to cause a drop in rank of $CT^0(q)$  or of its submatrices. Thus, a statement like: ``The network matrix $G^0(q)$ that is consistent with a given topology is generically identifiable from a given choice $\cal C$ of measurements'' will mean that $G^0(q)$ is identifiable for almost all choices of the  elements $G_{ij}^0(q)$ of $G^0(q)$ that are not known to be zero. 

We shall define this new notion of generic identifiability of the network in precise terms in Section~\ref{pathbased}. In order to give the reader an intuitive feeling for this notion, we illustrate it with the following example.\footnote{Starting in this example and for the remaining of the paper, we omit the dependence on $q$ whenever it creates
no confusion.}
 
\ni {\bf Example 1:} 
Consider a network whose topology is defined by the following network matrix:
$$G^0 = \left(\begin{array}{ccccc}0 & 0 & 0 & 0 & 0 \\G^0_{21} & 0 & 0 & 0 & 0 \\G^0_{31} & 0 & 0 & 0 & 0 \\0 & G^0_{42} & G^0_{43} & 0 & 0 \\0 & G^0_{52} & G^0_{53} & 0 & 0\end{array}\right) $$
and suppose we measure nodes 4 and 5 only. Simple calculations show that 
$$CT^0= \left(\begin{array}{ccccc}G^0_{42}G^0_{21} + G^0_{43}G^0_{31} & G^0_{42} & G^0_{43} & 1 & 0 \\
G^0_{52}G^0_{21} + G^0_{53}G^0_{31} & G^0_{52} & G^0_{53} & 0 & 1\end{array}\right)$$
Clearly, from $CT^0$ we can uniquely identify $G^0_{42}, G^0_{43}, G^0_{52}, G^0_{53}$. The remaining elements, $G_{21}$ and $G_{31}$ are then recovered from $CT^0$ by solving
$$\left(\begin{array}{cc}G^0_{42} & G^0_{43} \\G^0_{52} & G^0_{53}\end{array}\right)\left(\begin{array}{c} G_{21} \\ G_{31}\end{array}\right) = \left(\begin{array}{c}T^0_{41} \\ T^0_{51}\end{array}\right)$$
We conclude that $G^0$ is generically identifiable from measurements of nodes 4 and 5 only, because it is identifiable for almost all network matrices $G$ consistent with the topology, namely all  except those for which $G^0_{42} G^0_{53} = G^0_{52}G^0_{43}$, which is  a subset of measure zero.

\section{Motivating examples}\label{motivation}
In order to motivate the reader, we now analyze a few 3-node networks and show that the nodes that allow identification of the whole network depend entirely on the topology of the network, and that the whole network can often be identified from the measurements of a small subset of nodes. 

\begin{figure}[h]
 \centering
 \begin{tabular}{cc}
 \begin{minipage}{.4\textwidth}
\includegraphics[width=0.4\linewidth]{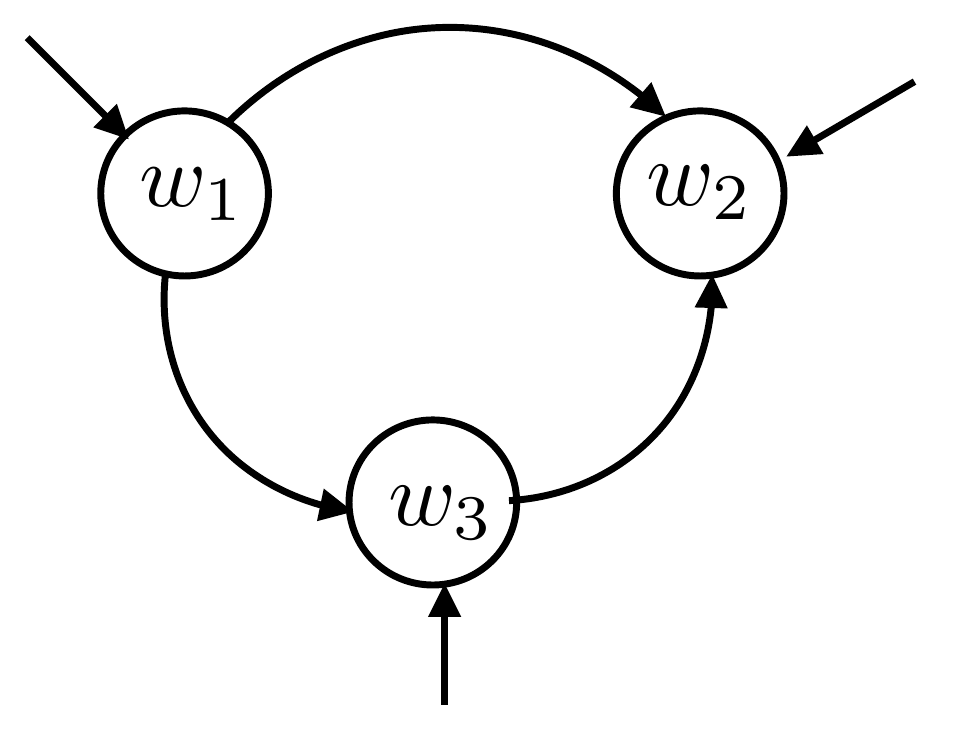} 
    \end{minipage}
  $\hspace{-40mm}  G^0 =\left[\begin{array}{ccc}0 & 0 & 0 \\ G_{21}^0 & 0 & G_{23}^0\\ G_{31}^0 & 0 & 0 \end{array}\right], $ \\     
$T^0 =  \left[\begin{array}{ccc}1 &0 & 0 \\ G_{21}^0 +  G_{23}^0G_{31}^0& 1 & G_{23}^0\\ G_{31}^0&0 &1 \end{array}\right] $
\end{tabular}
\caption{Example of network with three transfer functions where two nodes (2 and 3) need to be measured.}\label{fig:network1}
\end{figure}

Consider first a network with 3 unknown transfer functions represented in Figure \ref{fig:network1} and its corresponding true $G^0$ and true $T^0$. Calculations based on  (\ref{TzeroT2}) show that identification of all 3 transfer functions requires the measurement of nodes 2 AND 3, and that measuring node 1 yields no information.

\begin{figure}[h]
 \centering
 \begin{tabular}{cc}
   \begin{minipage}{.4\textwidth}
\includegraphics[width=0.4\linewidth]{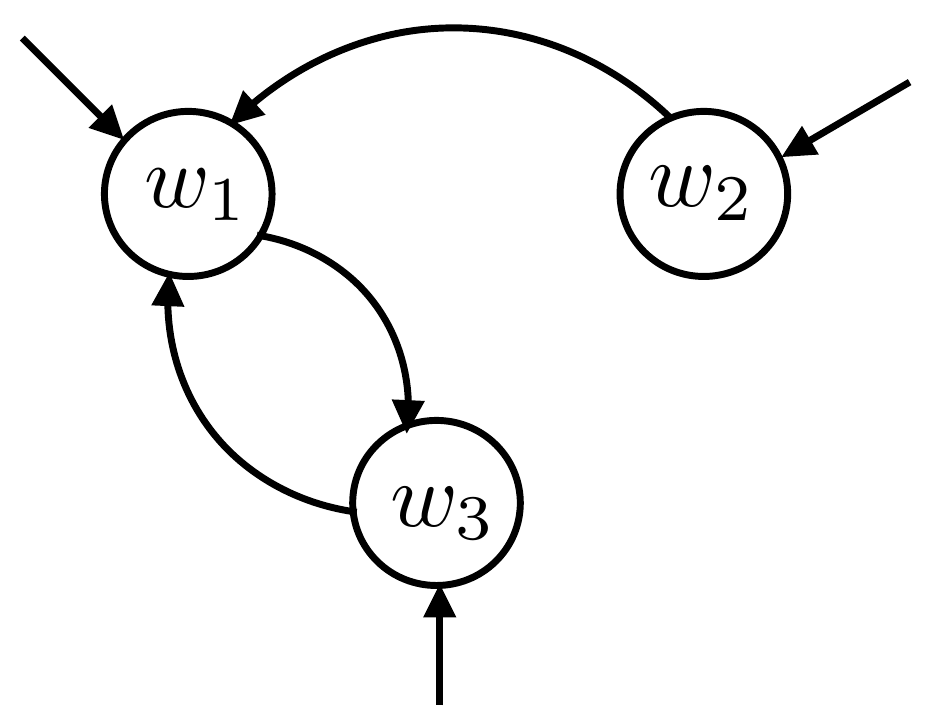} 
    \end{minipage}
 $\hspace{-40mm} G^0 =\left[\begin{array}{ccc}0 & G_{12}^0 & G_{13}^0 \\ 0 & 0 & 0 \\ G_{31}^0 & 0 & 0 \end{array}\right], $ \\
 $T^0 =  \frac{1}{\Delta} \left[\begin{array}{ccc}1 &G_{12}^0  & G_{13}^0 \\ 0 & 1- G_{13}^0G_{31}^0 & 0 \\ G_{31}^0  & G_{31}^0 G_{12}^0 & 1 \end{array}\right] $
\end{tabular}
\caption{Example of network with three transfer functions where measuring one node (1 or 3) is sufficient. We use $\Delta \bydef det(I-G^0)= 1 - G_{13}^0G_{31}^0$}\label{fig:network2}
\vspace{5mm}

 \centering
 \begin{tabular}{cc}
 \begin{minipage}{.4\textwidth}
\includegraphics[width=0.4\linewidth]{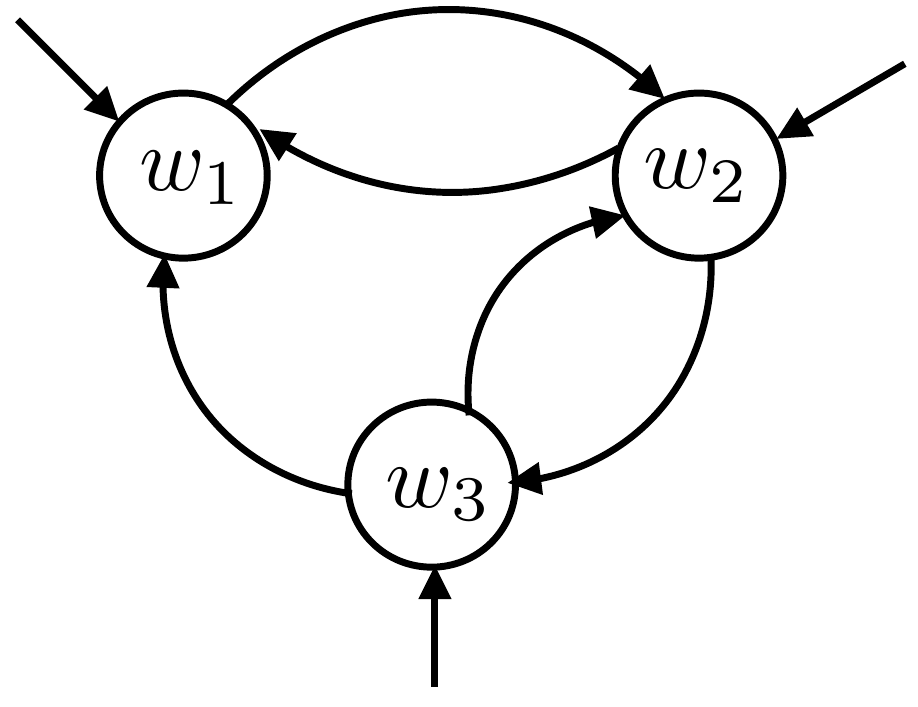} 
    \end{minipage}
 $\hspace{-40mm} G^0=\left[\begin{array}{ccc}0 & G_{12}^0 & G_{13}^0 \\ G_{21}^0 & 0 & G_{23}^0 \\ 0 & G_{32}^0 & 0 \end{array}\right], $ \\     
\hspace{-3mm} $T^0 \!\!=\!\!  \frac{1}{\Delta}\! \left[\begin{array}{ccc}1-G_{23}^0G_{32}^0  & G_{12}^0+G_{13}^0G_{32}^0    & G_{13}^0+ G_{12}^0G_{23}^0 \\ G_{21}^0 & 1 & G_{23} ^0+ G_{21}^0G_{13}^0 \\ G_{32}^0G_{21}^0  & G_{32}^0  & 1-G_{12}^0 G_{21}^0  \end{array}\!\! \right] $
\end{tabular}
\caption{Example of network with five transfer functions where measuring two nodes (1 and 2 or 1 and 3) is sufficient. We use $\Delta \bydef  1-G_{12}^0G_{21}^0 -  G_{23}^0G_{32}^0- G_{13}^0G_{21}^0G_{32}^0  $.}\label{fig:network3}
\end{figure}

By contrast, the identification of the 3 unknown transfer functions in the network represented in Figure \ref{fig:network2}  is possible by measuring just one node: node 1 OR node 3.

Finally,  in the network of Figure \ref{fig:network3}, all  5 transfer functions can be identified by measuring  just two nodes:  either nodes 1 AND 2, OR nodes 1 AND 3.  

These examples show  that the number and the choice of measurements that are necessary to identify the network depends not only on the number of unknown
transfer functions to be determined (the number of nonzero $G_{ij}$) but also on the topology of the network.

\section{Basic results}\label{basicresults}
Inspired by our analysis of 3-node networks, we now establish a number of basic results regarding the identifiability of general $L$-node networks from a reduced set of node measurements. In particular we show that 
measurements of some nodes is indispensable, and we  establish the minimum number of nodes that need to be measured for the identifiability of $G^0$. 

  We first introduce some  notations and we define some concepts from graph theory (see e.g. \cite{Diestel:00},  \cite{newman2010networks} ). Observe that the graph associated with  the network model is a graph with directed edges, i.e. if $G_{ij}^0$ is nonzero, it means that there is a directed edge from node $j$ to node $i$. Conversely, we say that a transfer matrix $G$ (or any matrix) is \emph{consistent} with a directed graph if $G_{ij} \neq 0$ only if there is an edge $(j,i)$.
  Observe that the presence of an edge does not require the corresponding entry $G_{ij}$ to be different from zero. \\
  {\bf Notations and definitions:} 
\begin{itemize}
\item $L$ = number of nodes;  
\item $p$ = number of measured nodes; 
\item $s$ = number of sinks, i.e. number of nodes with only incoming edges;  
\item $n$ = number of unknown transfer functions; 
\item $C$ = the $p \times L$ matrix that reflects the selection of nodes via $y(t) = C w(t)$: thus  each row of $C$ contains one element $1$ and $L-1$ elements $0$; 
\item $\cal C$ = the subset of nodes selected by $C$;
\item $G_{TS}^0$ = the restriction of the network matrix $G^0$ to the rows contained in a set $T$ and the columns contained in a set $S$;
\item $| \cal A|$ = cardinality of a set $\cal A$;
\item $N_i^+$ = set of out-neighbors  of node $i$, i.e. the set of nodes $j$ for which $G_{ji}^0\neq 0$;
\item  $d_i^+$ = $| N_i^+| $ = number of outgoing edges of node $i$; 
\item  a {\it walk}~denotes a series of adjacent directed edges (including trivial walks consisting of one node with no edge);
\item a {\it loop} is a walk whose terminal node coincides with the initial node;\footnote{Note that a loop is typically called {\it cycle} in graph theory.} 
\item  a {\it path} is a 
walk that never passes twice through the same node, i.e. a walk without loops. 
\item a directed graph is {\it weakly connected} if, for any partition of its vertices in two sets, there is at least one edge starting in one of the sets and ending in the other one.
\item a {\it tree} is a  graph that is weakly connected and has no loops even if one were to change the edges directions.
\end{itemize}

We can now establish the  following basic results. 
\begin{theorem}\label{basic1}
1) If $w_i$ is a source, then $G_{ij}^0 =0 \;\forall j$, $T_{ii}^0=1$ and $T_{ij}^0=0 \;\forall j\neq i$.
The measurement of a source does not add any linearly independent equation to the system of equations (\ref{TzeroT2}). The identification of a transfer function on an outgoing edge from a source $i$ requires that an external signal $r_i$ is applied at the source. \\
2) If $w_i$ is a sink, then $G_{ji}^0 =0 \;\forall j$, $T_{ii}^0=1$ and $T_{ji}^0=0 \;\forall j\neq i$. Identifiability of the network requires that all sinks be measured. The application of an external signal $r_i$ at a sink $i$ yields no information.\\
\ni{\bf Proof:}Ê 
1) The first part follows from the definition of a source and from the calculation of $T^0$ from such $G^0$ using 
(\ref{Tdef}). It then follows   that if $C$ selects a source, say $w_i$, then the corresponding equation of   (\ref{Tdef}) yields $T_{ii}^0=1$, which does not contribute any information for the identifiability of $G^0$. Finally, let $i$ be a source with an outgoing edge $G_{ki}^0$. It follows that $w_i=r_i + v_i$ and thus $w_k= G_{ki}^0(r_i +v_i) + x$, where $x$ contains only terms that do not involve $G_{ki}^0$. Hence the identification of $G_{ki}^0$ requires that $r_i \neq 0$.\\
2) The first part follows from the definition of a sink and from the calculation of $T^0$  using (\ref{Tdef}). Let node  $i$ be a sink and let node $k$ be connected to $i$ by a nonzero transfer function $G_{ik}^0$. Since node $i$ is a terminal node of the path from $k$ to $i$, no  node signal other than $w_i$ can give any information about  $G_{ik}^0$. On the other hand, applying an excitation signal $r_i$ to sink $i$ yields no information, since no path leaves node $i$.
\cqfd
\end{theorem}

We now make some observations concerning the number of useful equations that result from (\ref{TzeroT2})  for the computation of the $G_{ij}$. Each measured node contributes $L$ equations, but some of these may not yield any information, because they result in 1=1 or 0=0.

First we note that $L-1 \leq n \leq L(L-1)$, the first inequality being a consequence of the
connectedness of the graph. The number of equations is $p\times L$, so it is obvious
that we need $p \geq \frac{n}{L}$. It now follows from (\ref{TzeroT2}) and Theorem~\ref{basic1} that each  sink causes the appearance of one trivial equation $1=1$ in the sink's measurement, and also 
of one trivial equation $0=0$ at every other measurement. Hence the number of trivial equations 
caused by each sink equals $p$, and thus the total number of trivial equations
due to the existence of sinks is $ps$. Therefore  the number of useful equations is
at most $n_e = pL -  ps = p(L-s)$. We then have the following result.

\begin{theorem}\label{basic2}
Identifiability of the whole network requires measurement of all sinks plus at least $m$ more nodes such that 
\be \label{minnumber}
m+s \geq \frac{n}{L-s}
\ee
\ni{\bf Proof:} Given that the number of useful equations resulting from $p$ measurements is at most $p(L-s)$, identifiability of a  network with $n$ unknowns and $s$ sinks requires that $p(L-s)\geq n$, where $p=m+s$. This implies (\ref{minnumber}). \cqfd
\end{theorem}
\vspace{-3mm}

The next theorem yields a simple result for networks that have the structure of a tree.

\begin{theorem}\label{treetheo}
For a tree  
it is necessary and sufficient to measure all the sinks, assuming that none of the $G_{ij}$ that make up the tree are zero.\\
\ni{\bf Proof:} By Theorem~\ref{basic1}  it is necessary to measure all the sinks for any graph, so it remains to prove
sufficiency.  In a tree every sink will be the terminal node of a path.
Given that all transfer functions  $T_{ji}^0$ from any input $r_i$ to any sink is identifiable, in order to determine all the  $G_{kl}^0$ in
that path one can proceed backwards from the sink up to the root, since the transfer function from any given
$r_i$ to the sink is just the product of the $G_{kl}^0$ of each edge in the path from $r_i$ to the sink, none of which is zero by our assumption.  \cqfd \end{theorem}
\vspace{-6mm}
After a result for networks having a tree structure, the next result covers the case of loops. 
\begin{theorem}\label{theoloop}
Let the nodes $w_i$, $i \in {\cal I}$  form one  loop and assume that 
 no other loop in the graph contains any of these nodes. Suppose moreover that all the transfer functions involved in the loop are nonzero. Then measuring any one of these nodes is sufficient to identify all transfer functions in the loop. \\
\ni{\bf Proof:}Ê Let $\eta$ be the cardinality of $\cal I$ and consider,  without loss of generality,
that the nodes in the loop are labeled $i=1,\ldots, \eta$ sequentially, that is there is a link from each node
$i$ to node $i+1$, so that the $\eta$ transfer functions to be identified in the loop are
$G_{i+1,i}^0, \; i=1,\ldots , \eta-1$ and $G_{1,\eta}^0$. Since an external excitation signal is assumed to enter each node, input-output  identification provides all closed-loop transfer functions 
$T_{i,j}^0, i,j\in {\cal I}$, none of which are zero. Indeed,
\begin{eqnarray*}
&&\hspace{-5mm} T_{i,j}^0 = \frac{1}{\Delta} G_{i,i-1}G_{i-1,i-2} \ldots G_{1,\eta} G_{\eta,\eta-1} \ldots G_{j+1,j}, \; i < j  \\
&&\hspace{-5mm} T_{i,i}^0 = \frac{1}{\Delta} \\
&&\hspace{-5mm} T_{i,j}^0= \frac{1}{\Delta} G_{i,i-1}G_{i-1,i-2} \ldots  G_{j+1,j}, \; i > j
\end{eqnarray*}
where 
\begin{equation}\label{poiu}
\Delta = 1- G_{1,\eta} \Pi_{i=1,\ldots , \eta-1} G_{i+1,i} .
\end{equation}
Now, suppose we measure only the
``last" node $i=\eta$. Then  we have identified all the transfer functions $T_{\eta,j}^0$:
\begin{eqnarray}
&&\hspace{-5mm} T_{\eta,j}^0 = \frac{1}{\Delta} G_{\eta,\eta-1}G_{\eta-1,\eta-2} \ldots G_{j+1,j}, \; j=1,\ldots,\eta-1 \nonumber \\
&&\hspace{-5mm} T_{\eta,\eta}^0 = \frac{1}{\Delta} \label{uiop}
\end{eqnarray}
Now, notice that 
$$
G_{j+1,j} = \frac{T_{\eta,j}^0}{T_{\eta,j+1}^0}, \; j=1,\ldots,\eta-1
$$
which gives each one of the transfer functions in the path from node 1 to node $\eta$, that is
all transfer functions in the loop except $G_{1,\eta}$. Then this last transfer function can be obtained, from
(\ref{poiu}) and (\ref{uiop}), as 
$$
G_{1,\eta} = \frac{1}{\Pi_{i=1,\ldots , \eta-1} G_{i+1,i}} \prt{1- \frac{1}{T_{\eta,\eta}^0}}.
$$
The same reasoning holds if we measure  any other node, since it is just a question of relabeling the nodes. 
\cqfd
\end{theorem}
\vspace{-6mm}

\section{Path-based results}\label{pathbased}

In this section, we consider a specific node $i$ within the network and its out-going edges, i.e. the edges corresponding to the nonzero elements $G_{ji}^0$ within the network matrix. Recall that we denote by $N_i^+$ the corresponding set of out-neighbors of node $i$. We show that the generic identifiability of an edge\footnote{For reasons of brevity, we shall in future often refer to the {\it identifiability of an edge}, where this in fact means the identifiability of the transfer function corresponding to this edge.} or a group of edges leaving this node $i$ can be related to the structure of the paths from the corresponding out-neigbors  to the measured nodes.

Section \ref{sec:linalg} presents a linear algebraic reformulation of the identifiability problem, which involves submatrices of $T^0$. In Section \ref{sec:generic} we formally define the notion of generic identifiability, needed because of the risk of exceptional  rank drops in the submatrices of $T^0$.  Section \ref{sec:bac} establishes the link between the structure of paths in the network and the generic rank of certain submatrices of $T^0$. These relations are then used in Section \ref{sec:necsuf} to obtain necessary and sufficient conditions for identifiability  of  out-going edges of a specific node, and some corollaries are derived in Section \ref{sec:csq_nec_suf}.

\subsection{A linear algebraic reformulation}\label{sec:linalg}

Remember that $CT^0$ can be perfectly identified from $\{y,r\}$ data, and that therefore 
the transfer function $G_{ji}^0$ of an edge $(i,j)$ is identifiable if (\ref{TzeroT2}) implies $G_{ji} = G^0_{ji}$ for any $G$ \emph{consistent with the graph, i.e. with the topology}.
Define $\Delta \bydef G-G^0$, which is consistent with the graph if and only if $G$ is. The next Lemma shows how the identifiability of $G_{ji}^0$ depends on the kernel of a submatrix of the known $CT^0$, and hence on the rank of certain submatrices of $CT^0$.

\begin{lemma}\label{submatrixT0}
Let  $N^+_i$ be the set of out-neighbors of node $i$. Let $T^0_{C,N^+_i}$ denote the restriction of $T^0$ to the rows selected by $C$  and to the columns corresponding to $N^+_i$, and let $\Delta_{N^+_i,i} $ denote the restriction of the $i$-th column of $\Delta$ to the rows corresponding to $N^+_i$. Then $G_{ji}^0$ is identifiable from $CT^0$ if and only if
\begin{equation} 
\label{eq:reformulation_Delta_TCN}
T^0_{C,N^+_i} \Delta_{N^+_i,i} = 0 \Rightarrow \Delta_{ji} = 0.
\end{equation}
\end{lemma}
\ni {\bf Proof:} 
Substituting $G=G^0 + \Delta$ in (\ref{TzeroT2}) and remembering $T^0(I-G^0)=I$ shows that $G_{ji}$ is identifiable if and only if 
\be \label{eq:reformulation_Delta}
C T^0  \Delta  =0\Rightarrow \Delta_{ji}  = 0
\ee
for any $\Delta $ consistent with the graph. The left hand side of \eqref{eq:reformulation_Delta} actually consists of   $L$ independent linear systems of the form
$$
C T^0  \Delta_{:\ell}  = 0, \hspace{.5cm} \ell =1,\dots L. 
$$
The function $\Delta_{ji} $ only appears in one system, with $\ell=i$, and none of the functions appearing in that system appear in any other one. 
Hence $G_{ji}$  is identifiable if and only if
\be \label{eq:idcolumn}
C T^0  \Delta_{:i}  =0\Rightarrow \Delta_{ji}  = 0
\ee
for any $\Delta_{:i}$ consistent with the graph i.e. $\Delta_{ki} = 0$ if there is no edge $(i,k)\in G^0$. Remember that $G_{ki}$, and hence $\Delta_{ki}$, may be nonzero only if $k \in N^+_i$. 
We  use the notation $l \in \cal C$ to say that $l$ is a measured node. Condition \eqref{eq:idcolumn} can be rewritten as
\begin{equation}\label{eq:condition_columnn_detail}
\sum_{k\in N^+_i} T^0_{lk}\Delta_{ki} = 0, \forall l\in {\cal C} \Rightarrow \Delta_{ji} = 0,
\end{equation}
which is equivalent to (\ref{eq:reformulation_Delta_TCN}).
\cqfd
The identifiability of $G_{ji}^0$ is thus related to the rank of $T^0_{C,N^+_i}$ and, as will be seen, that of certain of its submatrices. We will see in Section \ref{sec:bac} how these are related to the topology.

\subsection{Generic properties}\label{sec:generic}

As seen in   Example 1, identifiability essentially depends on the known graph associated to $G^0$, except for network matrices $G^0$ that lie in subsets of measure $0$. We now formalize this notion, using an approach similar to 
that in \cite{Vanderwoude:96}. 
A \emph{rational transfer matrix parametrization consistent with a given graph} is defined in the following way.
For every edge $(j,i)$, set 
constants $p_{ij}, ~n_{ij}\in \mathbb{Z}_0^+ $, and parametrize $G_{ij}(z)$ by 
\begin{equation}\label{eq:def_param}
G_{ij}(z) = \kappa_{ij}z^{- p_{ij}}\frac{z^{ n_{ij}}+\beta_{ij}^{(n_{ij}-1)} z^{ n_{ij}-1}+\dots + \beta_{ij}^{(1)} z + \beta_{ij}^{(0)}  }{z^{n_{ij}}+\alpha_{ij}^{(n_{ij}-1)}  z^{n_{ij}-1}+\dots + \alpha_{ij}^{(1)} z + \alpha_{ij}^{(0)}}
\end{equation}
for real parameters $\kappa_{ij}$, $\alpha_{ij}^{m}$ and $\beta_{ij}^{m}$, ($1\leq m < n_{ij}$). For pairs $(j,i)$ not connected by an edge, let $G_{ij}(z) = 0$.  We collect all parameters $\kappa_{ij}$, $\alpha_{ij}^{m}, \beta_{ij}^{m}$ in a vector $P$, and denote by $G(P,z)$ the transfer matrix obtained by a specific parameter. 

We say that a property \emph{generically\footnote{The word \quotes{structurally} is also sometimes used, see e.g. \cite{lin1974structural}.} holds for a network matrix $G^0$} if, for any rational transfer matrix parametrization $G(P,z)$ consistent with the graph associated to $G^0$, the property holds for $G(P,z)$ for all parameters $P$ except possibly those lying on a zero measure set in $\Re^N$, where $N$ is the total number of parameters. \\

\ni {\bf notational remark.} \\
In the remainder of this paper, and in order to simplify notations, we will say that a property generically holds for $T^0= (I-G^0)^{-1}$ if for every parametrization $G(P,z)$ consistent with the graph associated to $G^0$, the property  holds for $T(P,z) := (I-G(P,z))^{-1}$ for all $P$ except possibly those lying on a zero measure set. We will use the same convention for properties holding for submatrices of $T^0$. \\

We have seen in Section \ref{sec:linalg} that identifiability is linked to the rank of certain matrices. Hence generic identifiability will be linked to the \emph{generic rank} of certain submatrices of $T^0$, i.e. the size of their largest generically nonsingular submatrix. This implies checking if the determinant of a matrix related to $G$ is generically nonzero. The following Lemma, when applied to $Q$ being the determinant of a matrix related to $G$, provides a convenient way of establishing this. See proof in Appendix \ref{appen:lemma_generic}.

\begin{lemma}\label{lem:structural}
Let $Q(.):\mathbb{C}^{L\times L}\to \mathbb{C}$ be an analytic function and consider a network matrix $G^0(z)$. 
If there exists a matrix $A\in \mathbb{C}^{L\times L}$  consistent with the graph associated to $G^0(z)$ such that $Q(A) \neq 0$, then $Q(G^0(z))$ is generically not identically zero as a function of $z$ (for polynomial or rational $Q(.)$, it  then has finitely many roots). Otherwise, $Q(G(z))\equiv 0$ for every $G(z)$ consistent with the graph. \cqfd
\end{lemma}

This leads to the following definition of a \emph{generically identifiable network matrix.}
\begin{definition} \label{identifiability}
A network matrix $G^0(z)$ is generically identifiable from a set of measured nodes defined by $C$ in (\ref{measures}) if, for any rational transfer matrix parametrization $G(P,z)$ consistent with the directed graph associated to $G^0(z)$, there holds

\beqna \label{genid}
C(I-G(P,z))^{-1} = C(I-\tilde G(z))^{-1} \Rightarrow G(P,z) = \tilde G(z)
\eeqna
for all parameters $P$ except possibly those lying on a zero measure set in $\Re^N$, where $\tilde G(z)$ is any network matrix consistent with the graph.
\cqfd
\end{definition}
This definition naturally extends to the generic identifiability of a specific transfer function in $G$ (or edge), or of a group of these.
We note that \eqref{genid} is exactly parallel to the definition of identifiability of a given transfer matrix, since \eqref{TzeroT1} can be rewritten $C(I-G^0)^{-1} = C(I-G)^{-1}$.

\subsection{Disconnecting sets, vertex-disjoint paths and matrix rank}\label{sec:bac}
\providecommand{\bac}{b_{\cal A\to \cal C}}

\begin{figure}
\centering
\includegraphics[scale = .5]{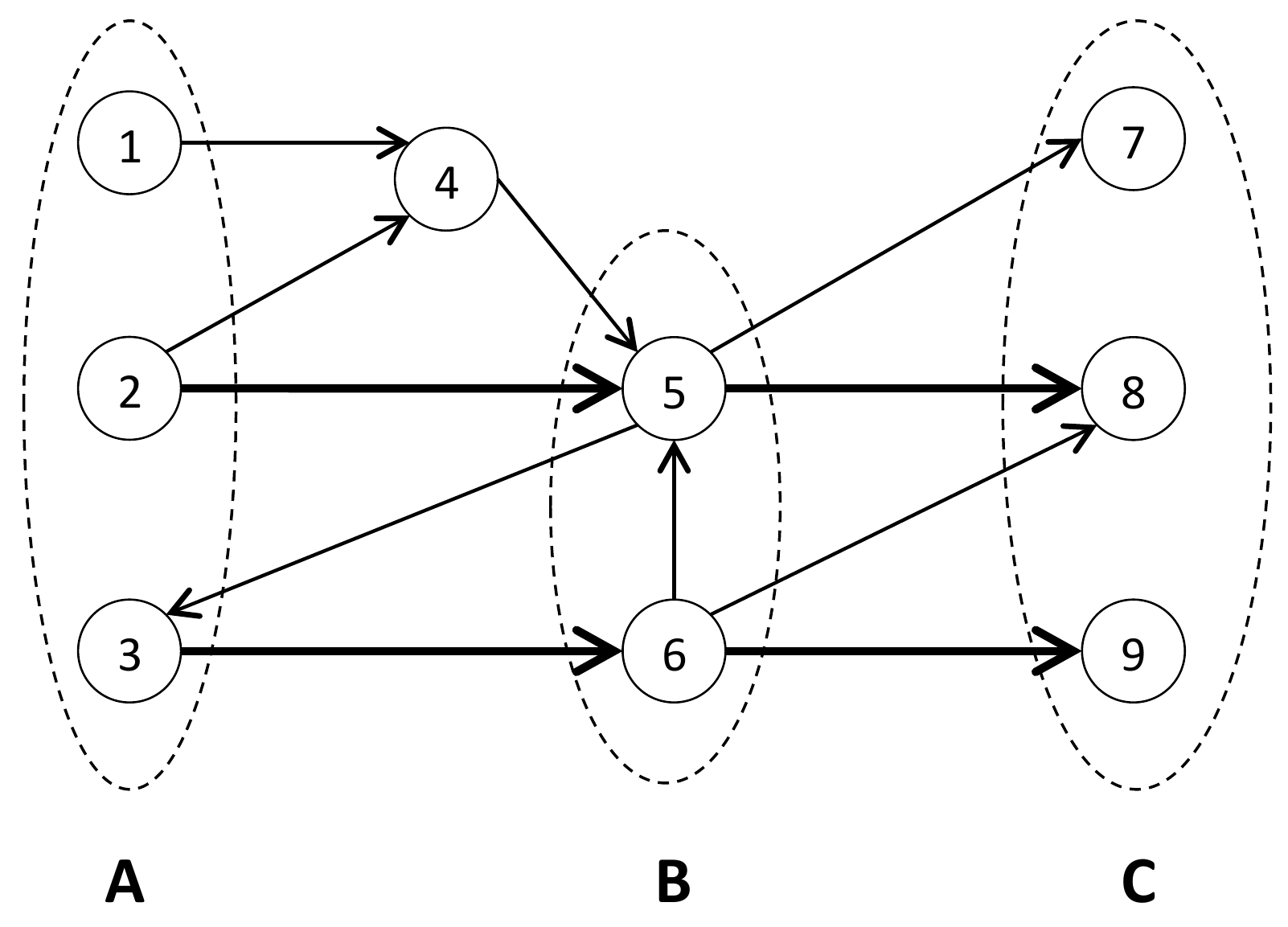}
\caption{Illustration of the notions of vertex-disjoint paths, disconnecting sets, and  Lemma \ref{lem:duality_bac}.
The highlighted paths $(2,5,8)$ and $(3,6,9)$ are vertex-disjoint paths from ${\cal A}=\{1,2,3\}$ to ${\cal C}=\{7,8,9\}$, and their existence implies $\bac \geq 2$. An alternative set of such paths is $(1,4,5,7) $ and $(3,6,9)$.
On the other hand, every path starting from ${\cal A}$ and arriving in ${\cal C}$ goes through ${\cal B}=\{5,6\}$, which is thus an $\cal A-\cal C$ disconnecting set. It follows then from Lemma \ref{lem:duality_bac} that $\bac \leq 2$ and thus that (i) $\cal B$ is a minimal $\cal A-\cal C$  disconnecting set, and (ii) there is no set of more than two vertex-disjoint paths from $\cal A$ to $\cal C$.}
\label{fig:duality_bac}
\end{figure}

We say that a group of paths are mutually vertex disjoint if no two paths of this group contain the same vertex. Consider two subsets  of nodes $\cal A$ and $\cal C$. We let $b_{\cal A\to \cal C}$ be the maximum number of mutually vertex disjoint paths starting in $\cal A$ and ending in $\cal C$.
We say that a set of nodes $\cal B$ is an \emph{$\cal A-\cal C$ disconnecting set} if every path starting in $\cal A$ and ending in $\cal C$ contains at least one node in $\cal B$, which implies that there would be no path from $\cal A$ to $\cal C$ if $\cal B$ were removed.   These notions are illustrated in Figure~\ref{fig:duality_bac}. Note that $\cal B$ can intersect $\cal A$ and/or $\cal C$. In particular, $\cal A$ and $\cal C$ are always  $\cal A-\cal C$ disconnecting sets.

The following Lemma, also illustrated in Figure \ref{fig:duality_bac}, links the  notions of disconnecting
sets and vertex disjoint paths.

\begin{lemma}\label{lem:duality_bac}
Consider two subsets  of nodes $\cal A$ and $\cal C$.  The maximum number $b_{\cal A\to \cal C}$ of mutually vertex disjoint paths from $\cal A$ to $\cal C$, is also the size of the smallest $\cal A-\cal C$ disconnecting set. Moreover, under the standing assumption that the network is weakly connected,
it can be computed in $O(Ln)$ operations.
\end{lemma}
\ni {\bf Proof:} 
The equality between $b_{\cal A\to \cal C}$ and the size of the smallest $\cal A-\cal C$
disconnecting set is the directed vertex disjoint version of  Menger's theorem, see e.g. \cite{Bohme&Frank&Harant:01}. 
Computing $b_{\cal A\to \cal C}$ can be recast as solving a max-flow problem, see for example Section 24.2 of \cite{Erickson:2014}. There exist many efficient ways of solving max-flow problems. Since the maximum flow is bounded by $L$, the classical Ford-Fulkerson Algorithm (see e.g. \cite[Section 10.5.1]{newman2010networks}), for example, terminates in $O(Ln)$ operations provided there are at least $L-1$ edges, which is the case if the network is weakly connected.\cqfd

The next lemma will be useful to use the bounds derived in terms of $b_{\cal A\to \cal C}$.

\begin{lemma}\label{lem:subaddititivy}
For any sets of nodes $\cal A,\cal A'$ and $\cal C$ there holds
$$ 
b_{\cal A\cup \cal A' \to \cal C} \leq b_{\cal A\to \cal C} + b_{\cal A'\to \cal C}
$$
\end{lemma}
\ni {\bf Proof:}
Consider a set of $b_{\cal A\cup \cal A' \to \cal C}$ vertex disjoint paths from $\cal A\cup \cal A'$ to $\cal C$. Let $b$ be the number among those starting from a node in $\cal A$. There must thus be at least $b_{\cal A\cup \cal A' \to \cal C} - b$ starting from $\cal A'$ (there can be more).
By definition, $b_{\cal A\to \cal C}\geq b$ because we have already found $b$ vertex disjoint paths from $\cal A$ to $\cal C$. Similarly, $b_{\cal A'\to \cal C}\geq b_{\cal A\cup \cal A' \to \cal C} - b$. 
Hence there holds
$$
  b_{\cal A\to \cal C} + b_{\cal A'\to \cal C} \geq b + b_{\cal A\cup \cal A' \to \cal C} - b = b_{\cal A\cup \cal A' \to \cal C}.
$$\cqfd

Our main result in this subsection is the establishment of the link between the generic rank of submatrices of  $T^0$  and the number of vertex-disjoint paths between the sets corresponding to the selected columns and rows of $T^0$. A  link between generic rank and vertex-disjoint paths was obtained in the pioneering paperÊ\cite{Vanderwoude:91}, where the  rank of the matrix $C(sI-A)^{-1}B$ of a system $\dot x = Ax +Bu$, $y = Cx$ was related with paths from inputs to outputs in a graph defined by $A,B$ and $C$. 

Our next Proposition differs from the main result of \cite{Vanderwoude:91} in several ways. First the paths defined in the graph associated to the matrix $C(sI-A)^{-1}B$ in \cite{Vanderwoude:91} are those of the whole network that connects the inputs of $C(sI-A)^{-1}B$ to its outputs, whereas we consider the paths connecting a subset of these inputs to a subset of the outputs. Secondly, the matrices $A, B, C$ appearing in 
$C(sI-A)^{-1}B$ are real matrices, while we examine the generic rank of a submatrix of $(I-G(P,z))^{-1}$ where $G(P,z)$ is a matrix of transfer functions. Finally, the definition of the nodes in the graph associated to $C(sI-A)^{-1}B$ differs from that used in this paper, because our nodes are linked by transfer functions; this means that if we were to represent $T^0(z) = (I-G^0(z)))^{-1}$ as a state space representation as is done in \cite{Vanderwoude:91}, then the nodes of the graph associated to our $G(z)$ would be a small subset of those associated to this state space representation. For all these reasons, the next Proposition is not just an application of the main theorem of \cite{Vanderwoude:91} and requires a specific proof.

\begin{proposition}\footnote{Remember the important notational convention adopted for $T^0$ and $T^0_{\cal{C},\cal{A}}(z).$} \label{prop:rank_cut}
Let $\cal A$, $\cal C$ be two sets of nodes of a directed graph associated to a network matrix $G^0(z)$. Let $T^0_{\cal{C},\cal{A}}(z)$ be the restriction of $T^0(z) = (I-G^0(z))^{-1}$  to the rows corresponding to $\cal C$ and columns corresponding to $\cal A$. Then the generic rank of  $T^0_{\cal{C},\cal{A}}(z)$ is $b_{\cal A\to \cal C}$.
\end{proposition}
\ni{\bf Proof:}
The proof will consist of two parts. The first one establishes that the rank is generically at least $\bac$, and uses the interpretation of $\bac$ in terms of the number of vertex-disjoint paths. \\
 \emph{Part 1: Generically $\Rank (T_{C,A}^0) \geq b_{\cal A\to \cal C}$ }\\
Select $\bac$ vertex-disjoint (directed) paths from $\cal A$ to $\cal C$, 
and let $A$ be the adjacency matrix of the directed graph consisting only of these paths, i.e. $A_{ij}=1$ if the edge $(j,i)$ is on one of the paths and $A_{ij}=0$ otherwise. It is then a standard result in graph-theory (see e.g. \cite[Section 6.10]{newman2010networks}) that $[A^k]_{ji}$ is the number of walks of length exactly $k$ from $j$ to $i$ in that graph. Since the graph consists of disjoint directed paths, this implies that (i) $A^{k}=0$ if $k$ is larger than the longest of the vertex-disjoint paths, and hence (ii)  $(I-A)^{-1}=\sum_{k=0}^\infty A^k$. As a result (iii) $[(1-A)^{-1}]_{ji}=\sum_{k=0}^\infty [A^k]_{ji}$ is the total number of walks of any length from $i$ to $j$ in the graph containing only the vertex disjoint paths. In particular, let now $\tilde {\cal A}\subseteq \cal A$ be the set of starting points of the paths, and $\tilde {\cal C} \subseteq \cal C$ the set of their arrival points, with obviously $|\tilde {\cal A}| = |\tilde {\cal C}| = \bac$. Therefore if $i\in \tilde {\cal A}$ and $j\in \tilde {\cal C}$ and if they are on the same path, then  $[(1-A)^{-1}]_{ji}=1$. Otherwise $[(1-A)^{-1}]_{ji}=0$ (as there is no walk from the origin of one path to the end of another one). The restriction $[(1-A)^{-1}]_{\tilde {\cal C}, \tilde {\cal A}}$ of $[(1-A)^{-1}]$ is thus a permutation matrix of size $\bac$, whose determinant is nonzero. By Lemma \ref{lem:structural} this implies that  $det({T_{\tilde C, \tilde A}^0}) $  is generically nonzero, implying that the rank 
of ${T_{\tilde C, \tilde A}^0}$ is generically $\bac$, and hence the generic rank of $T_{C,A}^0$ is at least $\bac$, since ${T_{\tilde C, \tilde A}^0}$ is a submatrix of $T_{C,A}^0$.

The proof of the second part  relies on the equivalent interpretation of $\bac$ in terms of the size of the minimal $\cal A-\cal C$ disconnecting set.\\
\emph{Part 2: Generically $\Rank (T_{C,A}^0) \leq b_{\cal A\to \cal C}$ }\\
Let $\cal B$ be an $\cal A-\cal C$ disconnecting set of minimal size $\bac$, the existence of which is guaranteed by Lemma \ref{lem:duality_bac}.
Let ${\cal S} \subset \{1,\dots,L\}$ be the set of nodes that can be reached by a path from a node in $\cal A$ without intersecting any node of the disconnecting set $\cal B$, and let ${\cal P} =  \{1,\dots,L\}\setminus(\cal S\cup \cal B)$. We have thus partitioned the $L$ nodes into 3 disjoint sets: $\cal P,\cal S$ and $\cal B$.
There holds $\cal A\subseteq  \cal S\cup \cal B$ (nodes in $\cal A$ are all in $\cal S$ except if they belong to $\cal B$). There also holds $\cal C\subseteq \cal P\cup \cal B$. Indeed, there would otherwise be a node of $\cal C$ in $\cal S$, meaning that it could be reached from a node in $\cal A$ without going through $\cal B$, in contradiction with $\cal B$ being an $\cal A-\cal C$ disconnecting set. 

After re-ordering of the indices, the matrices $G^0$ and $T^0$ can be rewritten as 
\beqna
G^0&=&\left(\begin{array}{ccc}
G_{PP}& G_{PB} & 0\\
G_{BP}& G_{BB} & G_{BS}\\
G_{SP}& G_{SB} & G_{SS}\\
\end{array}
\right)
\text{    and   }\\
T^0 &\bydef& (I-G^0)^{-1}=\left(\begin{array}{ccc}
T_{PP}& T_{PB} & T_{PS}\\
T_{BP}& T_{BB} & T_{BS}\\
T_{SP}& T_{SB} & T_{SS}\\
\end{array}
\right)
\eeqna
We focus on the rows $P$ and columns $S$ and $B$, keeping in mind that $T^0 = I+ G^0T^0$. There holds
\begin{align*}
T_{PS} = I_{PS} + [GT]_{PS} &= 0 +  G_{PP}T_{PS} + G_{PB} T_{BS} + 0 T_{SS}\\
T_{PB} = I_{PB} + [GT]_{PB} &= 0 +  G_{PP}T_{PB} + G_{PB} T_{BB} + 0 T_{SB},
\end{align*}
from which follows
\begin{equation}
\left(
\begin{array}{cc}
I-G_{PP} & 0 \\ 0 & I 
\end{array}
\right)
\left( 
\!\!\begin{array}{cc}
T_{PB} & T_{PS}\\
T_{BB} & T_{BS}  
\end{array}
\!\! \right)
=
\left( 
\!\! \begin{array}{c}
G_{PB}\\ I
\end{array}
\!\! \right)
\left( 
\begin{array}{cc}
T_{BB} & T_{BS}
\end{array}
\!\! \right).
\end{equation}

The right hand side of the equality has a rank at most $\bac$ because $(T_{BB}~~ T_{BS})$ has $\bac$ rows. The same holds thus true for the left-hand side. Observe now that the left-hand side is square and generically invertible; it is indeed invertible if we replace $G_{PP}$ by 0, and the generic invertibility then  follows  from Lemma \ref{lem:structural}. As a consequence, 
the rank of $T_{P\cup B,B\cup S}$, the second matrix of the left hand side, is also at most $\bac$. The claim of part 2 follows then from the fact that $T_{C,A}^0$ is a submatrix of $T_{P\cup B,B\cup S}$, because we have seen that $\cal A\subseteq \cal B\cup S$ and $\cal C\subseteq \cal P\cup \cal B$.
\cqfd
The result of Proposition~\ref{prop:rank_cut}
can intuitively be understood as follows. In the  system represented by $T^0$, an edge can carry a one-dimensional information about the effect of a given external excitation, and a vertex can only let a one-dimensional information about a given external excitation transit through it. 
Suppose first that the graph only consists of two paths starting in $\cal A$ and ending in $\cal C$. If the paths are vertex-disjoint, then nodes in $\cal A$ can transmit a two-dimensional information about a given external signal to those in $\cal C$, one dimension per path. On the other hand, if the two paths intersect in one vertex, only a one-dimensional information can transit through this vertex and reach $\cal C$. Proposition \ref{prop:rank_cut} extends this intuitive idea to graphs with more edges than just those on the paths and to larger number of paths. Since we know by Lemma \ref{lem:duality_bac} that the largest number of vertex-disjoint paths between two sets is the size of the smallest disconnecting set, this allows  characterizing exactly the dimension of the information transmitted, i.e. the rank of $T_{C,A}^0$. 
\subsection{Necessary and sufficient conditions for generic identifiability}
\label{sec:necsuf}

With the help of Proposition~\ref{prop:rank_cut} we can now derive one of the main results of this paper, namely necessary and sufficient conditions for the generic identifiability of transfer functions leaving a given node $i$. 
The reformulation of the identifiability of a transfer function leaving node $i$ by \eqref{eq:reformulation_Delta_TCN} naturally leads one to consider conditions for the generic identifiability of a group of edges leaving the same node $i$, as these are all related to the same matrix $T^0_{C,N^+_i}$.

\begin{theorem}\label{thm:nec_suf_group}
Let $N_i^*$ be a subset of $N_i^+$ and denote $\bar N_i^* \bydef N_i^+\setminus N_i^*$. The transfer functions corresponding to edges from $i$ to $N_i^*$ can generically all be identified when measuring nodes  $\cal C$ using the identified $CT^0$  if and only if the following two conditions hold:
\beqna \label{necsufni*}
 b_{N^*_i\to \cal C}\!\!&\!\!=\!\!&\!\!|N_i^*| \\
b_{N^+_i \to \cal C} \!\!&\!\!=\!\!&\!\! b_{N^*_i \to \cal C} + b_{\bar N_i^* \to \cal C} = |N_i^*|+ b_{\bar N_i^* \to \cal C}\label{eq:cond_b_nec_suf_group}
\eeqna
\end{theorem}
\vspace{2mm}

\providecommand{\nis}{N_i^*}
\providecommand{\bnis}{\bar N_i^*}
\ni {\bf Proof:}
Let us fix a $G$ consistent with the graph defined by $G^0$ and the corresponding $T=(I-G)^{-1}$. It follows from \eqref{eq:idcolumn} that we can recover the transfer functions of all edges from $i$ to $\nis$ if and only if the equality $CT \Delta_{:,i}=0$ implies $\Delta_{\nis,i}=0$ for every $\Delta_{:,i}$ for which $\Delta_{ki}=0$ for every $k\not\in N^+_i$. 
This can be rewritten as 
\begin{equation}\label{eq:sumTnisTbnis}
T_{C,\nis} \Delta_{\nis,i} + T_{C,\bnis}\Delta_{\bnis,i} =0 \rightarrow \Delta_{\nis,i}=0
\end{equation}
Observe first that $T_{C,\nis}$ must have rank $|\nis|$ for this condition to hold; otherwise one could find a $\Delta_{\nis,i}\neq 0$ for which $T_{C,\nis} \Delta_{\nis,i}=0$, which, with $\Delta_{\bnis,i}=0$, would contradict the condition. 
We can rewrite \eqref{eq:sumTnisTbnis} as 
$$
T_{C,\nis} \Delta_{\nis,i} = - T_{C,\bnis}\Delta_{\bnis,i}  \rightarrow \Delta_{\nis,i}=0
$$
If the image sets of $T_{C,\nis}$ and  $T_{C,\bnis}$ have a nontrivial intersection, then we could find $\tilde \Delta_{\nis,i}\neq 0$ and $\tilde \Delta_{\bnis,i}\neq 0$ such that 
$T_{C,\nis} \tilde \Delta_{\nis,i} = -  T_{C,\bnis} \tilde \Delta_{\bnis,i}\neq  0$, and the condition is not satisfied. On the other hand, if the image sets of $T_{C,\nis}$ and  $T_{C,\bnis}$ have no nontrivial intersection, then  $T_{C,\nis} \Delta_{\nis,i} = - T_{C,\bnis}\Delta_{\bnis,i}$ implies both  $T_{C,\nis} \Delta_{\nis,i}=0$ and  $T_{C,\bnis} \Delta_{\bnis,i}=0$. When $T_{C,\nis}$ has rank $|\nis|$, the former equality implies $ \Delta_{\nis,i}=0$.

We have thus shown that the transfer functions of edges from $i$ to $\nis$ can all be identified from $CT$  if and only if (i) $\Rank(T_{C,\nis}) = |\nis|$ and (ii) the image sets of $T_{C,\nis}$ and $T_{C,\bnis}$ have no nontrivial intersection, i.e. are linearly independent. The latter condition is equivalent to $\Rank(T_{C,(\nis\cup \bnis)})= \Rank(T_{C,\nis}) + \Rank (T_{C,\bnis})$. 
For any $G$ (and corresponding $T$) consistent with the graph associated to $G^0$,  this equality, together with   $\Rank (T_{C,\nis})= |\nis|$, are thus necessary and sufficient  for the identifiability of the transfer functions corresponding to the edges from  $i$ to $\bar N_i^*$.  This is in particular the case for any matrix $G(P,z)$ in any parametrization of the transfer matrices consistent with that graph. Generic identifiability of the edges from $i$ to $\bar N_i^*$ is thus equivalent to $\Rank(T^0_{C,\nis})= |\nis|$ and  $\Rank(T^0_{C,(\nis\cup \bnis)})= \Rank(T^0_{C,\nis}) + \Rank (T^0_{C,\bnis})$ holding generically, and the equivalence with \eqref{necsufni*} and \eqref{eq:cond_b_nec_suf_group}  then follows from  Proposition \ref{prop:rank_cut}.
\cqfd

\ni   {\bf Comment:} Condition  \eqref{eq:cond_b_nec_suf_group} can also be formulated as
\be \label{ineq:cond_b_nec_suf_group} 
b_{N^+_i \to \cal C} \geq b_{N^*_i\to \cal C} + b_{\bar N_i^* \to \cal C}  = |N_i^*|+ b_{\bar N_i^* \to \cal C},
\ee
because it follows  from the sub-additivity Lemma \ref{lem:subaddititivy} and $N^+_i= N^*_i\cup \bar N^*_i$ that 
$$
b_{N^+_i \to \cal C} \leq b_{N^*_i\to \cal C} + b_{\bar N_i^* \to \cal C}.
$$
This formulation will be used in the proof of Corollary~\ref{cor:suf_partial}.

An immediate corollary of Theorem~\ref{thm:nec_suf_group} is obtained when one considers all out-neighbors of node $i$.

\begin{corollary}\label{cor:necandsuf}
The transfer functions  from node $i$ to its out-neighbors $N_i^+$ can generically all be identified from $CT^0$ if and only if $ b_{N^+_i\to \cal C}= d_i^+ \bydef |N_i^+| $.
\end{corollary}
\ni {\bf Proof:}
The result follows from Theorem \ref{thm:nec_suf_group} applied to $N_i^*=N_i^+$, in which case $\bar N_i^*$ is an empty set. 
\cqfd 
Corollary \ref{cor:necandsuf} can be intuitively understood in the following way. We want to recover $d_i^+$ transfer functions of edges leaving $i$, so we need a $d_{i}^+$-dimensional information about the effect of $r_i$. Moreover, the information we have comes from the out-neighbors of $i$ and arrives at our measured nodes $\cal C$. Hence the recovery will be possible if and only if a $d_i^+$-dimensional information is transmitted from these out-neighbors to $\cal C$, which requires  $b_{N^+_i\to \cal C}= d_i^+ $ by Proposition \ref{prop:rank_cut}. In case we only want to recover the transfer function of the edges arriving at a subset $N_i^*$ of the out-neighbors of $i$ as in Theorem \ref{thm:nec_suf_group}, then the situation is more complex because the information received at $\cal C$ from $N_i^*$ is mixed with information about other edges leaving $i$. One  then has to check if the specific information about $N_i^*$ can be isolated in all the information arriving at $\cal C$ from the out-neighbors of $i$, which is what condition \eqref{eq:cond_b_nec_suf_group} is about. It can indeed be interpreted as requiring all the loss of information-dimension from $N_i^+$ to $\cal C$ to concern exclusively information about $N_i^+\setminus N_i^*$, leaving that about $N_i^*$ intact.

\begin{figure}[h]
 \centering
\includegraphics[width=0.7\linewidth]{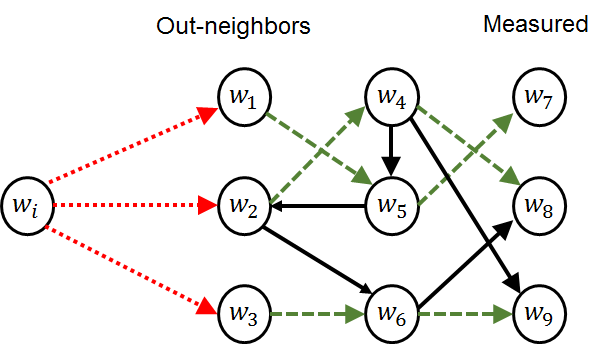}     
\caption{Example illustrating Corollary~\ref{cor:necandsuf2}: 3 vertex-disjoint dashed green paths to the 3 measured nodes; the 3 dotted red edges are identifiable.}
\label{fig:nscillustration}
\end{figure}

Given the definition of $ b_{\cal A \to \cal C}$, an alternative formulation of the previous result is as follows.
\begin{corollary}\label{cor:necandsuf2}
The transfer functions  from node $i$ to  its out-neighbors $N_i^+$ can generically all be identified from $CT^0$  if and only if there exist vertex-disjoint directed paths\footnote{The vertex-disjoint condition applies also for the departure and arrival nodes.} leaving all out-neighbors of $i$ and arriving at the measured nodes defined by $C$. \cqfd
\end{corollary}
\vspace{-3mm}
Corollary~\ref{cor:necandsuf2} is illustrated by the  example in Figure~\ref{fig:nscillustration}. Remember that known external signals $r_i$ are applied to each node, which we have not added on the figure for visibility reasons. Node $i$  has three outgoing nodes, each of which has a vertex-disjoint directed path to the measured nodes 7, 8 and 9, namely the paths $(1,5,7), (2,4,8)$ and $(3,6,9)$; they are represented by dashed green arrows. As a result, the dotted red transfer functions $G_{1i}^0, G_{2i}^0$ and $G_{3i}^0$ can all be identified from these three measured nodes.

We stress that the sufficient condition in Corollary \ref{cor:necandsuf2} 
does not require all paths from the nodes $N_i^+$ to the measured nodes to be disjoint, but only the existence of a set of mutually disjoint paths. In other words, there may very well exist many other paths than those used in the condition, and there is no requirement on those, nor on their intersections with those used in the condition.   For example,  Figure~\ref{fig:nscillustration} illustrates that the conditions of Corollary~\ref{cor:necandsuf2} apply even though  node 2 has another path to node 8, namely $(2,6,8)$ which has a common node with the path $(3,6,9)$.

Particularizing Theorem \ref{thm:nec_suf_group} to a set $N_i^*$ consisting of a single node immediately leads to a necessary and sufficient condition for identifying a single transfer function.

\begin{theorem}\label{thm:nec_suf_indiv}
Consider an edge $(i,j)$ and its corresponding transfer function $G_{ji}^0$, and let $N^+_i$ be the set of out-neighbors of $i$. The transfer function $G_{ji}^0$ can be generically uniquely identified
by measuring the nodes $\cal C$ if and only if 
\be \label{theo52cond}
b_{N^+_i \to \cal C} = b_{N^+_i\setminus \{ j \} \to \cal C} +1 
\ee
\end{theorem}

\ni {\bf Proof:}
The result follows directly from Theorem \ref{thm:nec_suf_group} applied to $N_i^*=\{j\}$, taking into account the fact that $b_{j\to \cal C}$ is 1 if there is a path from $j$ to $\cal C$ and 0 otherwise.
\cqfd
\subsection{Additional results}\label{sec:csq_nec_suf}
In this Subsection we present several results that apply to specific cases and
that can be directly derived from the previous results. We start by giving 
conditions for identifying a group of edges that are only sufficient
(not necessary) but that 
are simpler than the ones given in Theorem \ref{thm:nec_suf_group}.

\begin{corollary}\label{cor:suf_partial}
Consider a node $i$, and let $N^*_i\subseteq N^+_i$ be a subset of its out-neighbors with $|N^*_i|=d^*_i$. Suppose in addition that the two following conditions hold

(i) There exist $d^*_i$ vertex disjoint directed paths joining the nodes of $N^*_i$ to the   measured nodes   $\cal C$,

(ii) There is no path from any node of $N^+_i\setminus N^*_i$ to any node of $\cal C$.

\ni Then all transfer functions  from node $i$ to nodes in $N_i^*$ can be generically identified from the measured nodes.
\end{corollary}
\ni {\bf Proof:}
Since $|N_i^*|=d^*_i$, there holds $b_{N_i^* \to \cal C}\leq d^*_i$. Hence it follows from condition (i) that $b_{N_i^* \to \cal C}= d^*_i$.
Condition (ii) implies that $b_{N_i^+\setminus N^*_i \to \cal C} = 0$. Now $b_{N_i^+\to \cal C} $ is by definition always larger than or equal to $b_{N_i^*\to \cal C}$ because $N_i^*\subseteq N_i^+$. It also follows from Lemma \ref{lem:subaddititivy} that 
\be \label{eq:cor53bis}
b_{N^+_i \to \cal C} \leq b_{N^*_i\to \cal C}  +b_{N^+_i\setminus N_i^* \to \cal C}= b_{N^*_i\to \cal C}
\ee
where the last equality follows from $b_{N_i^+\setminus N^*_i \to \cal C} = 0$. Combining this with $b_{N_i^+\to \cal C}\geq b_{N^*_i\to \cal C}$ yields the desired result by Theorem  \ref{thm:nec_suf_group}. This implies that 
\be \label{eq:cor53}
b_{N^+_i \to \cal C} \geq b_{N^*_i\to \cal C} = b_{N^*_i\to \cal C} +b_{N^+_i\setminus N_i^+ \to \cal C}.
\ee
The result then follows  from Theorem \ref{thm:nec_suf_group} and Lemma \ref{lem:subaddititivy}.
\cqfd

The next two results concern the generic identification of the whole network, starting with
a rather simple but very telling necessary condition.  It is  related again to the need of obtaining information of sufficiently high dimension for the out-neighbors of every node.

\begin{corollary}\label{cor:highestdegree}
The network can be generically identified from $CT^0$ only if  $\cal C$ contains at least as many nodes as the highest out-degree present in the network.
\end{corollary}
\ni {\bf Proof:}
This follows from a direct application of the condition of Corollary~\ref{cor:necandsuf2} to a node $i$ with the highest out-degree.
\cqfd

Our next (and final) result in this subsection deals with the number of nodes
that are necessary and sufficient for identification of the whole network: it
shows that we never need to measure all $L$ nodes to secure network identifiability. 
It also confirms that, without any knowledge of the topology, we need to measure at least all but one of the nodes, as we cannot exclude the possibility of the graph being fully connected. 
To prove it, we first need the  following Lemma.

\begin{lemma}\label{allmin1}
Suppose we measure the set ${\cal C}_k$ defined as containing all nodes except $k$. If the network cannot be generically fully identified from $C_kT^0$, then there exists a node $k'$ with $N^+_{k'} \supseteq  N^+_k \cup \{k\}$, and thus $d^+_{k'}>d^+_k$.\\ 
\ni{\bf Proof:} 
If a node $i$ does not have $k$ as out-neighbor, then all its outgoing edges are generically identifiable. Indeed, all out-neighbors belong to ${\cal C}_k$ and are thus all connected to ${\cal C}_k$  by trivial zero-length vertex disjoint paths. Suppose now that node $i$ has $k$ as out-neighbor. If $k$ has an out-neighbor $j$ that is not an out-neighbor of $i$, then there exists $d^+_i$ vertex disjoint paths from $N_i^+$ to ${\cal C}_k$: the path $(k,j)$, and $d^+_i-1$ zero-length trivial paths from the other out-neighbors of $i$ to themselves (since they belong to ${\cal C}_k$). Hence the network is generically fully identifiable. So if the network is not generically fully identifiable, then $k$ and all its out-neighbors must be out-neighbors of $i$, which proves the claim with $k'=i$. \cqfd
\end{lemma}

\begin{theorem}\label{Lminus1}
$p=L-1$ is sufficient for identifiability, in the sense that there always exists a set of $L-1$ measured nodes allowing to generically fully identify the network. In particular, measuring all nodes except one of those with the highest out-degree is always sufficient.
In a fully connected network (that is, $n=L(L-1)$), $p=L-1$  is also necessary.
\\
\ni{\bf Proof:}
Necessity for the fully connected case is obvious,
since to identify $n=L(L-1)$ unknowns we need at least $p=\frac{n}{L}=L-1$ measurements.
To prove sufficiency, consider the set ${\cal C}_k$ defined in Lemma~\ref{allmin1}, where $k$ is a node with the maximal out-degree. It then follows from  Lemma~\ref{allmin1}  that this ${\cal C}_k$ allows generically full identification of the network, for otherwise there would be a node $k'$ with a higher out-degree, which is a contradiction.
\cqfd
\end{theorem}
In this Section, we have started our identifiability analysis by looking at a given node and its outgoing edges. We have given necessary and sufficient conditions for the identifiability of one specific outgoing edge, or a subset of outgoing edges, or all of them. These conditions are based on the existence of disjoint paths from these outgoing edges to the measured nodes. In particular, our results are useful to decide which nodes need to be measured if one wants to identify a particular transfer function: see Theorem~\ref{thm:nec_suf_indiv}. In addition, we have shown that it is never necessary to measure all $L$ nodes of a network, but that $L-1$ measures are sufficient. Several of our results and examples have actually shown that special structures within the network often allow one to identify the network using a much smaller number of measurements than $L-1$: see   e.g. Theorem~\ref{theoloop} and Example~\ref{fig:nscillustration}.

\section{Measurement-based results}\label{measurementbased}

\begin{figure}[h]
 \centering
\includegraphics[width=0.7\linewidth]{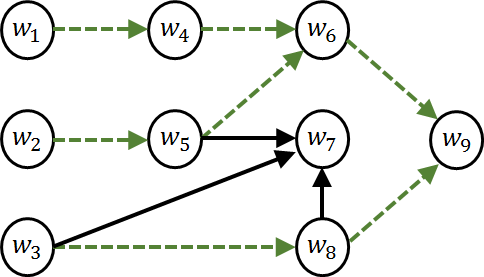}     
\caption{Example illustrating Theorem~\ref{measurednode}: all transfer functions on the dashed green edges can be generically identified from the measurement of node 9.}
\label{fig:migillustration}
\end{figure}

In this section, instead of starting from a given  node and its outgoing edges, we look at the converse approach. We consider a measured node, or a set of measured nodes, and we examine which transfer functions are identifiable from that measured node or from this set of measured nodes. In the first result we consider a single measured node.
\begin{theorem}\label{measurednode}
Let $j$ be a measured node, and consider a node $i$ that has a  path to node $j$. Then all transfer functions along that path can generically be identified if there is no other walk that connects $i$ to $j$.
\end{theorem}
\ni {\bf Proof:}
Let  $N^*_i$ of Theorem~\ref{thm:nec_suf_group} contain only the out-neighbor of node $i$ that is on the path to $j$ mentioned in the theorem, and let $\cal C$ contain only $j$. By the assumption in 
the statement, there is no path from any node in $N^+_i\setminus N^*_i$ to $j$ since this would constitute another path from $i$ to $j$. The result then follows by applying Theorem~\ref{thm:nec_suf_group} to the successive nodes along the path from $i$ to $j$. 
\cqfd
Theorem~\ref{measurednode} is illustrated by the  example in Figure~\ref{fig:migillustration}; remember again that known signals $r_i$ are added to each node, which are not shown on the figure. It follows from this theorem that the  7 transfer functions on the dashed green-colored  paths can all be generically identified from the measurement of node 9. If in addition node 7 is also measured, then the 10 transfer functions of the network can all be generically identified from the two measured nodes 7 and 9.
 The intuition behind Theorem \ref{measurednode} is that for each edge of the path we need to recover a specific one-dimensional information about the effect of the input at its starting node. The presence of the path from $i$ to $j$ guarantees that this one-dimensional information reaches $j$, and the absence of another walk to $j$ guarantees that it is not mixed with other information about the same input, i.e. about other edges leaving the node.
\\

The following result extends Theorem~\ref{measurednode} by providing a necessary and sufficient condition for the generic identifiability of all transfer functions on a path to a single measured node. 

\begin{theorem}\label{measurednodef}
Let $j$ be the only measured node and consider a node $i$ that has a path to node $j$ - let's call it path $\cal P$.
All transfer functions along $\cal P$ can be (generically) identified if and only if any other walk from $i$ to $j$ contains $\cal P$ as a prefix.\footnote{A path ${\cal P}_1$ is a prefix to another path 
${\cal P}_2$ if the initial nodes of ${\cal P}_2$ are those of ${\cal P}_1$.}
\end{theorem}

\providecommand{\lp}{\ell_{\cal P}}
\providecommand{\lw}{\ell_{\cal W}}
\ni {\bf Proof:} 
\emph{Necessity:} Suppose there is a walk $\cal W$ from $i$ to $j$ that does not contain $\cal P$ as prefix. Since they both start from $i$, $\cal P$ and $\cal W$ begin by a common part, possibly reduced to node $i$ without any edge.   Let $k$ be the last node of this initial common part, that is, $\cal P$ and $\cal W$ are identical until $k$ and different afterwards. This last common  node $k$ cannot be $j$ for otherwise $\cal P$ would be a prefix of $\cal W$. 
Hence there is a node after $k$ along $\cal P$, which we call $\ell_{\cal P}$, and a node after $k$ along $\cal W$, which we call we call $\ell_{\cal W}$.
We
apply Theorem \ref{thm:nec_suf_indiv} to the edge $(k,\ell_{\cal P})$. Clearly $b_{N_k^+\to \cal C}=1$ because ${\cal C}=\{j\}$ 
contains only one node. Moreover, there is by definition a walk from $\ell_{\cal W}$ to $j$ and,
since $\ell_{\cal W}\in  N^+_k$, there is a path from 
$\ell_{\cal W}\in  N^+_k\setminus \{\lp \}$ to $j$, 
so that $b_{N^+_k\setminus\{\lp \} \to \cal C}=1$. It follows then from Theorem \ref{thm:nec_suf_indiv} that  $(k,\ell_{\cal P})$ cannot be generically identified.

\emph{Sufficiency:} Suppose now there exists a node $k\neq j$ on $\cal P$ and its successor $\lp$ is such that the edge $(k,\lp)$ is not generically identifiable. Clearly, $b_{\lp\to \cal C}=1$. Hence it follows from Theorem \ref{thm:nec_suf_indiv} that $b_{N^+_k\setminus\{\lp \} \to \cal C}=1$, which means there exists another 
neighbor, that we call $\lw$, from which there is a path $\cal P'$ to ${\cal C}=\{j\}$. We can then build a walk $\cal W$ from $i$ to $j$ by aggregating (i) the 
restriction of $\cal P$ to its first nodes until it arrives at $k$, (ii) the edge $(k,\lw)$ and (iii) the path $\cal P'$ from $\lw$ to $j$, and this walk does not contain $\cal P$ as a prefix.
\cqfd
Finally, the results of Section~\ref{pathbased} allow us to produce  a necessary and sufficient condition for the generic identifiability of all edges of the network from a given set of measured nodes, i.e. a given choice of $C$. This is another main result of this paper.

\begin{theorem}\label{thm:nec_suf_global}
All edges of a network can generically be identified    if and only if $b_{N^+_i\to \cal C} = d^+_i$ for every $i$. 
\end{theorem}
\ni {\bf Proof:} The result follows immediately from Corollary \ref{cor:necandsuf}, or the equivalent Corollary \ref{cor:necandsuf2}, applied to all nodes.\cqfd
Theorem \ref{thm:nec_suf_global} can be put in other (more intuitive) words as
follows: all edges can be identified if and only if for every node $i$ there exist
$d^+_i$ vertex-disjoint paths from the set of neighbors of $i$ to  the nodes of $\cal C$.

\section{Algorithmic Complexity}\label{algocom}

Our results allow determining whether a given set $\cal C$ of measured node allows recovering a specific edge $(i,j)$, a specific set of edges,
or all edges in the network. Let us now analyze the algorithmic complexity of these issues. We have seen in Lemma \ref{lem:duality_bac} that $\bac$ can be computed in $O(Ln)$ for any sets $\cal A,\cal C$ using, for example, the Ford-Fulkerson algorithm.

It follows from Theorem \ref{thm:nec_suf_global} that checking if all edges can be identified can be achieved by computing $b_{N^+_i\to \cal C}$ for the $L$ nodes $i$, at a cost $L.O(Ln)= O(L^2n)$. 
If we only want to determine if a specific edge $(i,j)$ can be identified, then by 
Theorem \ref{thm:nec_suf_indiv} we can achieve this by comparing $b_{N^+_i \to \cal C}$ with $b_{N^+_i\setminus \{ j \} \to \cal C}$, the computation of which has a cost $O(Ln)$. Finally, suppose we are given a $\cal C$ and we want to determine the exact set of edges that can be identified. We then need to compute $b_{N^+_i \to \cal C}$ for each of the $L$ nodes $i$ and $b_{N^+_i\setminus \{ j \} \to \cal C}$ for each of the $n$ edges $(i,j)$, at a total cost of $(L+n)O(nL)=O((L+n)nL=O(n^2L)$ if we assume that the network is weakly connected, so that $L\leq n+1$.

\section{Conclusions}\label{conclusion}

The results so far on the global identifiability of a network of dynamical systems have been built on the assumption that all nodes are measured. In this paper, we have addressed the network identifiability problem in the situation where not all nodes are measured, but where they are all excited by a known external excitation signal. We have first shown that network identifiability with partial node measurements is impossible without knowledge about the topology. We have then developed an identifiability theory for a network matrix that is based on the topology of its associated graph, and not on the particular numbers that appear in the unknown network matrix. This has led us to define and exploit the notion of generic identifiability of a network matrix.

We have first shown that the node measurements needed for network identifiability depend entirely on the topology of the network. In doing so, we have observed that the measurement of all sinks are indispensable. 

We have then provided a series of results on identifiability. Some of these are based on looking at a particular node and its out-neighbours, and their paths to measured nodes; others have addressed the question of which transfer functions can be identified from the measurement of a particular node or a subset of nodes.

Our first main result, based on  the first approach, is a necessary and sufficient condition for identifiability of one edge, a set of  edges, or all edges leaving a particular node. Our second main result is a necessary and sufficient condition for identifiability of all transfer functions of the network from a selected set of measured nodes. We have also shown that these necessary and sufficient conditions
can be checked by algorithms that run in polynomial time, an important feature for large networks.

An interesting outcome of our work is that networks can often be identified by measuring only a small subset of  nodes.

Future research questions will include the search for a reduced  set of measured nodes that allow identification of the whole network, as well as the search for informative experiment designs.



\appendix

\ni {\bf Proof of Lemma \ref{lem:structural}} 
\label{appen:lemma_generic}

We first show that the absence of $A$ such that $Q(A) \neq 0$ implies $Q(G(z))\equiv 0$ for any $G$ consistent with the graph associated to $G^0(z)$. Indeed,  if there is a $G(z)$ such that $Q(G(z))\not \equiv 0$, then there is a $z^*$ such that $Q(G(z^*))\neq 0$ and we obtain the desired $A$ by taking $G(z^*)$.

We now show that the existence of $A$ such that $Q(A) \neq 0$ implies $Q(G^0(z))\neq 0$ generically. Consider a parametrization $G(P,z)$ of rational transfer functions consistent with the graph associated to $G^0$, and let  $\tilde Q:(P,z)\to \tilde Q(P,z) = Q(G(P,z))$ as a function of both $z$ and the parameters collected in $P$. Suppose, to obtain a contradiction, that the implication does not hold, that is there exists a nonzero-measure set $\mathbb{P}_0$ of parameters $P$ such that $Q(G(P,z))$, as a function of $z$, is identically zero. This implies that $\tilde Q(P,z) = 0$ for every couple $(P,z)\in  \mathbb{P}_0\times \mathbb{C}$, a set whose measure is also nonzero. 
Now, it follows from the assumption on $Q$ and the parametrization by rational functions that $\tilde Q$ is analytic. And it is a classical result that analytic functions that are not identically zero vanish only  on a zero-measure set \cite{Mityagin:15}. 
In particular, the fact that $\tilde Q$ vanishes on $\mathbb{P}_0\times \mathbb{C} $ implies that it is identically 0. 

We now show that this contradicts the existence of $A$ consistent with the graph for which $Q(A) \neq 0$.  Observe indeed that $A = G(P^*,1)$ for a parametrization $P^*$ defined by letting $k^*_{ij}= A_{ij}$, and $\beta^{t*}_{ij}=\alpha^{t*}_{ij}=0$ for every $i,j$ and $t$. Hence $\tilde Q(P^*,1) = Q(G(P^*,1))= Q(A) \neq 0$, in contradiction with $\tilde Q \equiv 0$. Therefore, the existence of $A$ implies that $Q(G(P,z))$  is identically zero (as a function of $z$)  only on a zero-measure set of parameters $P$. The last part of the result follows from the fact that single-variable polynomials have a finite number of roots when they are not identically zero, and the same holds for rational functions. \cqfd

\newpage
%
\begin{IEEEbiography}[{\includegraphics[width=1in,height=1.25in,clip,keepaspectratio]{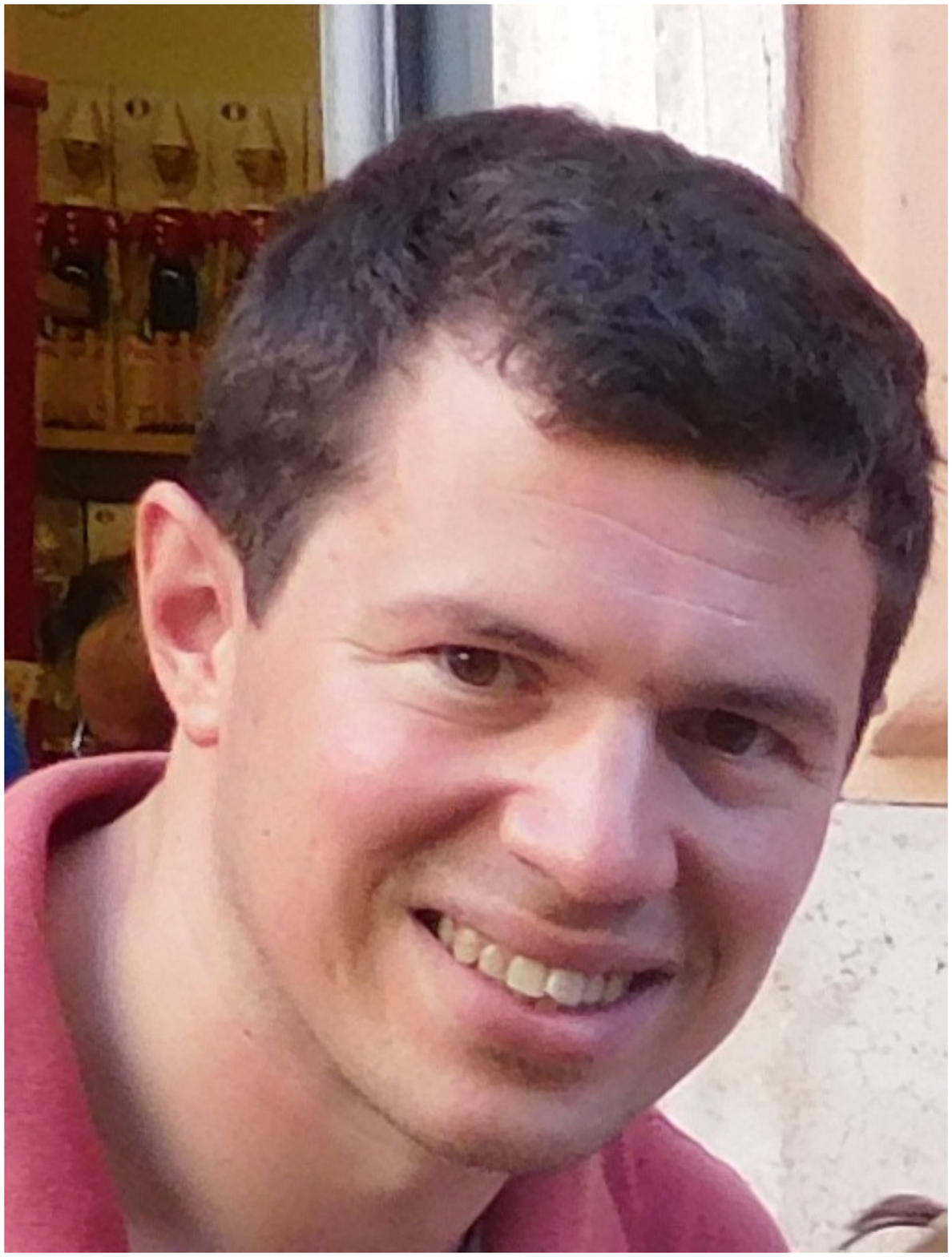}}]{Julien M. Hendrickx}
Julien M. Hendrickx received an engineering degree in applied mathematics and a PhD in mathematical engineering from the Universit\'e catholique de Louvain, Belgium, in 2004 and 2008, respectively.

He has been a visiting researcher at the University of Illinois at Urbana Champaign in 2003-2004, at the National ICT Australia in 2005 and 2006, and at the Massachusetts Institute of Technology in 2006 and 2008. He was a postdoctoral fellow at the Laboratory for Information and Decision Systems of the Massachusetts Institute of Technology 2009 and 2010, holding postdoctoral fellowships of the F.R.S.-FNRS (Fund for Scientific Research) and of the Belgian American Education Foundation. Since September 2010, he is a faculty member of the Universit\'e catholique de Louvain, in the Ecole Polytechnique de Louvain.

Doctor Hendrickx is the recipient of the 2008 EECI award for the best PhD thesis in Europe in the field of Embedded and Networked Control Systems, and of the Alcatel-Lucent-Bell 2009 award for a PhD thesis on original new concepts or application in the domain of information or communication technologies.
\end{IEEEbiography}
\newpage

\begin{IEEEbiography}[{\includegraphics[width=1in,height=1.25in,clip,keepaspectratio]{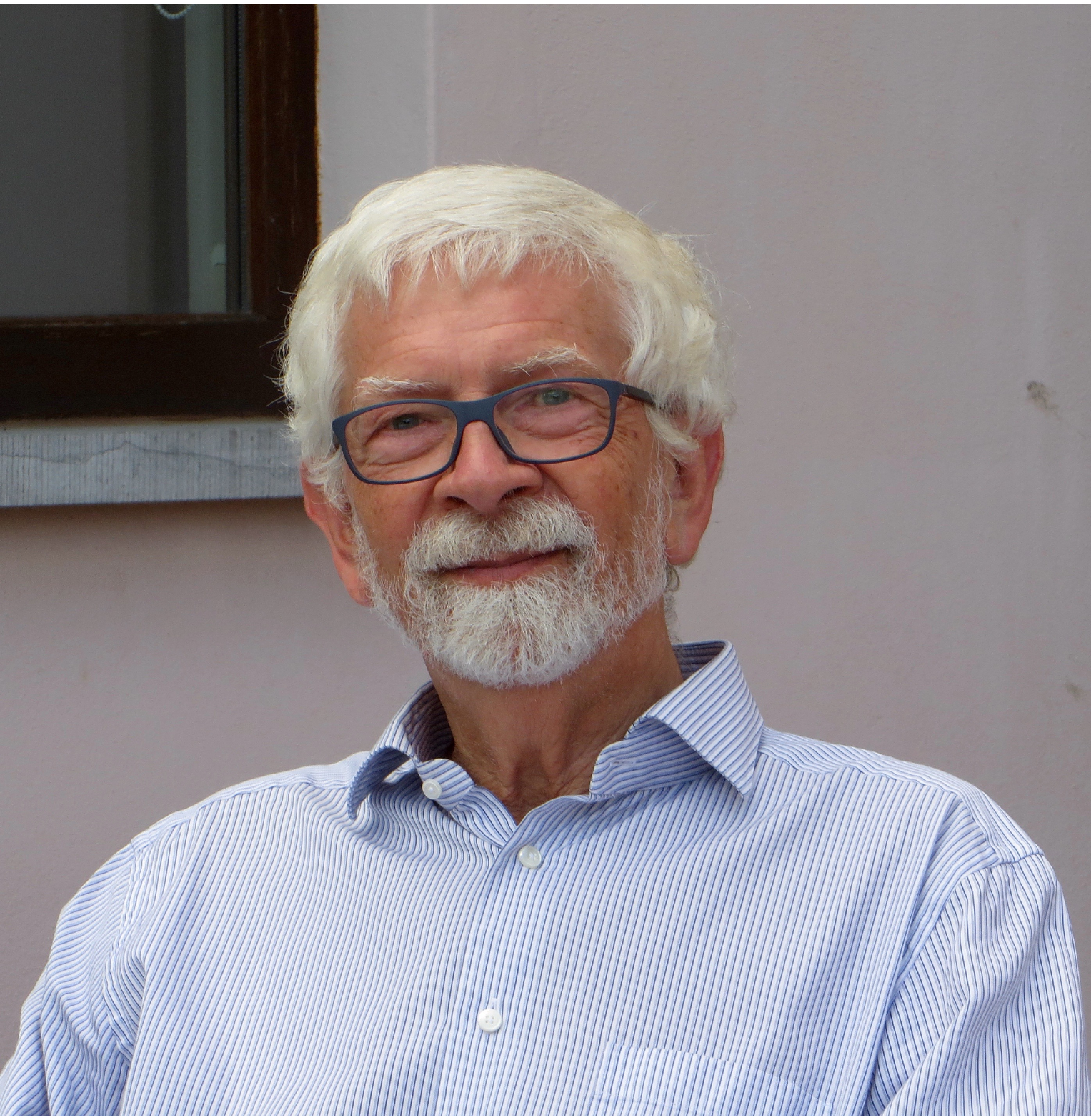}}]{Michel Gevers}
Michel Gevers obtained an Electrical Engineering degree from the University of Louvain, Belgium, in 1968, and a Ph.D. degree from Stanford University, California, in 1972, under the supervision of Tom Kailath. He is an IFAC Fellow, a Fellow of the IEEE, a Distinguished Member of the IEEE Control Systems Society. He holds  Doctor Honoris Causa degrees from the University of Brussels and from Link\"oping University, Sweden. He has been President of the European Union Control Association (EUCA) from 1997 to 1999, and Vice President of the IEEE Control Systems Society in 2000 and 2001. 

Michel Gevers is Professor Emeritus at the Department of Mathematical Engineering of the University of Louvain, Louvain la Neuve, and scientific collaborator at the Free University Brussels (VUB). He has been for 20 years the coordinator of the Belgian Interuniversity Network DYSCO (Dynamical Systems, Control, and Optimization) funded by the Federal Ministry of Science.  He has spent long-term visits at the University of Newcastle, Australia, and the Technical University of Vienna, and was a Senior Research Fellow at the Australian National University from 1983 to 1986. His present research interests are in system identification, experiment design for identification of linear and nonlinear systems, and identifiability and informativity issues in networks of linear systems. 

Michel Gevers has been Associate Editor of Automatica, of the IEEE Transactions on Automatic Control, of Mathematics of Control, Signals, and Systems (MCSS) and Associate Editor at Large of the European Journal of Control. He has published about 280 papers and conference papers, and two books: "Adaptive Optimal Control - The Thinking Man's GPC", by R.R. Bitmead, M. Gevers  and V. Wertz (Prentice Hall, 1990), and  "Parametrizations in Control, Estimation and Filtering Problems: Accuracy Aspects", by M. Gevers and G. Li  (Springer-Verlag, 1993). 
\end{IEEEbiography}

\begin{IEEEbiography}[{\includegraphics[width=1in,height=1.25in,clip,keepaspectratio]{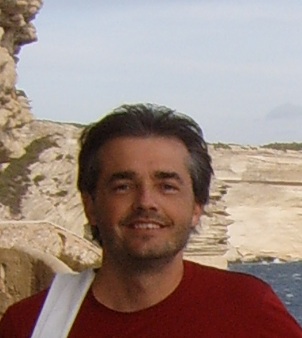}}]{Alexandre Bazanella}
Alex Bazanella received his PhD degree in Electrical Engineering in 1997.
He is currently a Full Professor with the Department of Automation
and Energy of Universidade Federal do Rio Grande do Sul, in Porto Alegre, Brazil.
His main research interests are presently in system identification
and data-driven control design, but he has also authored a number of papers
in nonlinear systems theory, particularly its application to electric
power systems and electrical machines.
He is the author of two books:
"Control Systems - Fundamentals and Design Methods" (in portuguese),
and "Data-driven Controller Design: the H2 Approach" (Springer, 2011).
Dr. Bazanella has served as associate editor of the IEEE Transactions on
Control Systems Technology from 2002 to 2008 and as an Editor of
the Journal of Control, Automation and Electrical Systems from 2008 to 2012.
He has held visiting professor appointments
at Universidade Federal da Para'ba in 2001 and at 
Universit\'e catholique de Louvain, where he spent a sabbatical year in 2006.
Dr. Bazanella is a Senior Member of the IEEE, and is currently vice-president
of the Brazilian Automation Society and vice-coordinator of the
area of Electrical Engineering for CAPES, with the Brazilian Ministry
of Education.
\end{IEEEbiography}


%



\end{document}